\documentclass[12pt]{article}
\usepackage{graphicx, amsmath, mathrsfs}
 \usepackage{color}
\usepackage{amsfonts}
\usepackage{verbatim}
\usepackage{latexsym}

\makeatletter
\@addtoreset{equation}{section}

\begin{document}

\newcommand{\E}{\mathbb{E}}
\newcommand{\PP}{\mathbb{P}}
\newcommand{\RR}{\mathbb{R}}
\newcommand{\LL}{\mathbb{L}}

\newtheorem{theorem}{Theorem}[section]
\newtheorem{remark}{Remark}[section]
\newtheorem{lemma}{Lemma}[section]
\newtheorem{coro}{Corollary}[section]
\newtheorem{defn}{Definition}[section]
\newtheorem{assp}{Assumption}[section]
\newtheorem{expl}{Example}[section]
\newtheorem{prop}{Proposition}[section]

\newcommand\tq{{\scriptstyle{3\over 4 }\scriptstyle}}
\newcommand\qua{{\scriptstyle{1\over 4 }\scriptstyle}}
\newcommand\hf{{\textstyle{1\over 2 }\displaystyle}}
\newcommand\hhf{{\scriptstyle{1\over 2 }\scriptstyle}}

\newcommand{\eproof}{\indent\vrule height6pt width4pt depth1pt\hfil\par\medbreak}

\def\a{\alpha} \def\g{\gamma}
\def\e{\varepsilon} \def\z{\zeta} \def\y{\eta} \def\o{\theta}
\def\vo{\vartheta} \def\k{\kappa} \def\l{\lambda} \def\m{\mu} \def\n{\nu}
\def\x{\xi}  \def\r{\rho} \def\s{\sigma}
\def\p{\phi} \def\f{\varphi}   \def\w{\omega}
\def\q{\surd} \def\i{\bot} \def\h{\forall} \def\j{\emptyset}

\def\be{\beta} \def\de{\delta} \def\up{\upsilon} \def\eq{\equiv}
\def\ve{\vee} \def\we{\wedge}
\def\vv{\varepsilon}

\def\F{{\cal F}}
\def\T{\tau} \def\G{\Gamma}  \def\D{\Delta} \def\O{\Theta} \def\L{\Lambda}
\def\X{\Xi} \def\S{\Sigma} \def\W{\Omega}
\def\M{\partial} \def\N{\nabla} \def\Ex{\exists} \def\K{\times}
\def\V{\bigvee} \def\U{\bigwedge}

\def\1{\oslash} \def\2{\oplus} \def\3{\otimes} \def\4{\ominus}
\def\5{\circ} \def\6{\odot} \def\7{\backslash} \def\8{\infty}
\def\9{\bigcap} \def\0{\bigcup} \def\+{\pm} \def\-{\mp}
\def\la{\langle} \def\ra{\rangle}

\def\tl{\tilde}
\def\trace{\hbox{\rm trace}}
\def\diag{\hbox{\rm diag}}
\def\for{\quad\hbox{for }}
\def\refer{\hangindent=0.3in\hangafter=1}

\newcommand\wD{\widehat{\D}}
\newcommand{\ka}{\kappa_{10}}

\title{
\bf Multilevel Monte Carlo theta EM scheme for SDDEs with small noise\thanks{Supported by NSFC(Nos., 11561027, 11661039), NSF of Jiangxi(Nos., 20181BAB201005, 2018ACB21001).}
 }

\author{
{\bf  Li Tan$^{a,b}$\, and \, Chenggui Yuan$^c$}\\
 \footnotesize{$^{a}$ School of Statistics, Jiangxi University of Finance and Economics, Nanchang, Jiangxi, 330013, P. R. China}\\
 \footnotesize{$^{b}$ Research Center of Applied Statistics, Jiangxi University of Finance and Economics,}\\
  \footnotesize{ Nanchang, Jiangxi, 330013, P. R. China}\\
\footnotesize{$^c$ Department of Mathematics, Swansea University, Swansea, SA2 8PP, U. K. }\\
\footnotesize{Email: tltanli@126.com, C.Yuan@swansea.ac.uk}}

\date{}

\maketitle

\begin{abstract}
In this paper, a multilevel Monte Carlo theta EM scheme is provided for stochastic differential delay equations with small noise. Under a global Lipschitz condition, the variance of two coupled paths is derived. Then, the global Lipschitz condition is replaced by one-sided Lipschitz condition, in order to guarantee the moment finiteness of numerical scheme, a modified multilevel Monte Carlo theta EM scheme is put forward and the second moment of two coupled paths is estimated.


\end{abstract}

\section{Introduction}
Small noise stochastic differential equations (SDEs) are widely used in economics, finance, computational fluid dynamics, ecology, population dynamics and etc, and many customized numerical methods have been developed for small noise SDEs with the aim of improving efficiency \cite{ahs15, mt97, rw06}. The mostly used numerical methods are the Euler-Maruyama (EM) scheme, the Monte Carlo method, the Milstein method and the Runge-Kutta method. There are a lot of results for numerical schemes of  SDEs under the global Lipschitz condition, see \cite{kp92, pb11}, etc. Since the global Lipschitz condition is too strong for most equations, more and more works on SDEs with the non-global Lipschitz conditions are established in recent years. For SDEs under the non-global Lipschitz condition, the numerical schemes may not reproduce the behaviour of exact solutions \cite{hmy07}, or the moments of numerical solutions may even explode in a finite time\cite{hjk11}. Thus, the classical numerical schemes are modified or improved to guarantee the finiteness of numerical solutions or to improve efficiency under non-global Lipschitz conditions, for example the tamed EM scheme \cite{hjk11}, the truncated EM scheme \cite{m15}, the theta EM scheme \cite{kps91, ty18}, the tamed Milstein method \cite{wg13}, the multilevel Monte Carlo method \cite{g081, g082, hjk13, glm17}.

In \cite{ahs15}, the authors proposed a multilevel Monte Carlo EM method for stochastic differential equations with small noise, analyzed the variance between two coupled paths, and discovered that the computational complexity of multilevel Monte Carlo method combined with standard EM scheme was lower than the standard Monte Carlo. However, the results of \cite{ahs15} are obtained under the global Lipschitz condition. If the global Lipschitz condition weakens to one-sided Lipschitz condition, will it remain the same property? Motivated by \cite{ahs15}, we combine multilevel Monte Carlo method with the theta EM scheme and consider the variance between two coupled paths for stochastic differential delay equations (SDDEs) with small noise under global Lipschitz condition. Then, we replace the global Lipschitz condition by the one-sided Lipschitz condition, give a modified multilevel Monte Carlo theta EM scheme in order to guarantee the moment finiteness of the scheme. The second moment of two coupled paths is estimated under one-sided Lipschitz condition.

 Throughout this paper, unless otherwise specified, we let $(\Omega,
  {\cal{F}},\{{\cal{F}}_{t}\}_{t\ge 0}, \mathbb{P})$ be a complete probability space
  with a filtration
  $\{{\cal F}_t\}_{t\ge 0}$ satisfying
 the usual conditions.
  %(i.e. it is increasing and right continuous while
  %${\cal{F}}_0$ contains all
 %$\mathbb{P}$-null sets).
 Let  $W(t)=$ $(W_1(t),\ldots,W_d(t))^T$ be an $d$-dimensional Brownian motion defined on
 the probability space. Let $\tau>0$ be a delay. Consider the stochastic differential delay equation with small noise of the form
\begin{equation}\label{2.0}
{\mbox d}X^{\vv}(t)=f(X^{\vv}(t),X^{\vv}(t- \tau)){\mbox d}t+ \vv g(X^{\vv}(t),X^{\vv}(t- \tau)){\mbox d}W(t),t \geq 0
\end{equation}
with initial data $X(\theta)=\xi(\theta), \theta\in[-\tau,0]$, where $\vv\in(0,1)$ and
\begin{displaymath}
 f:\RR^a \times \RR^a
\rightarrow \RR^a  \mbox{ and }g: \RR^a \times \RR^a \rightarrow \RR^{a \times d}.
\end{displaymath}
In the following, we will analyze multilevel Monte Carlo EM solution of \eqref{2.0} under the global Lipschitz condition and one-sided Lipschitz condition respectively.

\section{SDDEs with Global Lipschitz Condition}\label{sec2}
  We shall impose the following hypothesis:

 \begin{description}
\item[(H)] Both $f$ and $g$ satisfy the global Lipschitz condition.  That is, there exists an $\alpha>1$ such that
\begin{displaymath}
|f(x,y)-f (\bar{x},\bar{y})| + |g(x,y)-g (\bar{x},\bar{y})| \le \alpha(|x-\bar{x}|+|y-\bar{y}|)
\end{displaymath}
for all $x, y, \bar x, \bar y \in \RR^a$. Moreover, for all $x, y\in \RR^a$
\begin{displaymath}
|\nabla f(x,y)|^2\vee|\nabla^2 f(x,y)|^2\le \alpha.
\end{displaymath}
\end{description}

\begin{lemma}
Let assumption (H) hold. Then, for any $T>0$ and $p\ge 2,$ we have
$$
\E \left[\sup_{0\le t\le T}|X^\vv(t)|^p\right] \le C.
$$
\end{lemma}

\begin{remark}\label{onem}
{\rm Assumption (H) implies the existence and uniqueness of equation \eqref{2.0}. Moreover, if  (H) holds, then for any $x, y \in \RR^a$
\begin{displaymath}
|f(x,y)| + |g(x,y)|\le \beta(1+|x| +|y|)
\end{displaymath}
where $\beta=\max\{\alpha,|f(0,0)|,|g(0,0)|\}$, and for any $x, y, \bar x, \bar y \in \RR^a$
\begin{displaymath}
\langle x-\bar{x},f(x,y)-f(\bar{x},\bar{y})\rangle \le \bar{\alpha}(|x-\bar{x}|^2+|y-\bar{y}|^2)
\end{displaymath}
where $\bar{\alpha}=\frac{1}{2}+\alpha^2$.
}
\end{remark}

\subsection{The theta EM Scheme}\label{sec2.1}
We now introduce theta EM scheme for \eqref{2.0}. Given any time $T>0$, assume that  there exist two positive integers such that $h=\frac{\tau}{m}=\frac{T}{M}$, where $h\in (0,1)$ is the step size. For $n=-m, \cdots, 0$, set $X_h^\vv(t_n)=\xi(nh)$; For $n=0, 1, \cdots,M-1$, we form
\begin{equation}\label{discrete}
\begin{split}
X_h^\vv(t_{n+1})=&X_h^\vv(t_n)+\theta f(X_h^\vv(t_{n+1}), X_h^\vv(t_{n+1-m}))h\\
&+(1-\theta) f(X_h^\vv(t_{n}), X_h^\vv(t_{n-m}))h+\vv g(X_h^\vv(t_{n}), X_h^\vv(t_{n-m}))\Delta W(t_n),
\end{split}
\end{equation}
where $t_n=nh$, $\Delta W(t_n)=W(t_{n+1})-W(t_n)$. Here $\theta\in [0,1]$ is an additional parameter that allows us to control the implicitness of the numerical scheme. For $\theta=0$, the theta EM scheme reduces to the EM scheme, and for $\theta=1$, it is exactly the backward EM scheme. For a given $X_h^\vv(t_n)$, in order to guarantee a unique solution $X_h^\vv(t_{n+1})$ to \eqref{discrete}, the step size is required to satisfy $\theta h<\frac{1}{\bar{\alpha}}$ according to the monotone operator \cite{zei}, where $\bar{\alpha}$ is defined as in Remark \ref{onem}. In addition, to guarantee the moment finiteness of numerical solutions, we also require $h\theta<\frac{1}{6\beta}$ in this section. Thus, in Section \ref{sec2}, we set $h^*\in\left(0,\frac{1}{\theta(\bar{\alpha}\vee 6\beta)}\right)$, and let $h, h_l\in(0,h^*]$ for $\theta\in(0,1]$, while for $\theta=0$, we only need  $h, h_l\in(0,1)$, where $h_l$ is a step size defined in Section \ref{sec2.2}.

We find it is convenient to work with a continuous form of a numerical method. Rewrite \eqref{discrete} with a continuous form as follows:
\begin{equation}\label{continuous}
\begin{split}
&X_h^\vv(t)-\theta f(X_h^\vv(t),X_h^\vv(t-\tau))h=\xi(0)-\theta f(\xi(0),\xi(-\tau))h\\
&+\int_0^tf(X_h^{\vv}(\eta_h(s)), X_h^{\vv}(\eta_h(s-\tau))){\mbox d}s+\vv\int_0^tg(X_h^{\vv}(\eta_h(s)), X_h^{\vv}(\eta_h(s-\tau))){\mbox d}W(s),
\end{split}
\end{equation}
where $\eta_h(s)=\left\lfloor s/h\right\rfloor h$.

\begin{lemma}\label{0pmoment}
Let assumption (H) hold. Then, for any $T>0$ and $p\ge 2,$ we have
$$
\E \left[\sup_{0\le t\le T}|X_{h}^\vv(t)|^p\right] \le C.
$$
\end{lemma}
{\bf Proof.} Denote $Y_h^\vv(t)=X_h^\vv(t)-\theta f(X_h^\vv(t),X_h^\vv(t-\tau))h$. For any $t\in[0,T]$, by \eqref{continuous} and the Burkholder-Davis-Gundy (BDG) inequality, we get
\begin{equation}\label{july}
\begin{split}
\E \left[\sup_{0\le s\le t}|Y_{h}^\vv(s)|^p\right]\le& C|Y_{h}^\vv(0)|^p+Ct^{p-1}\E\int_0^t|f(X^{\vv}_{h}(\eta_{h}(s)), X^{\vv}_{h}(\eta_{h}(s-\tau)))|^p{\mbox d}s\\
&+C\vv^p\mathbb{E}\left[\sup_{0\le s\le t}\left|\int_0^sg(X^{\vv}_{h}(\eta_{h}(u)), X^{\vv}_{h}(\eta_{h}(u-\tau))){\mbox d}W(u)\right|^p\right]\\
\le& C+Ct^{p-1}\E\int_0^t|f(X^{\vv}_{h}(\eta_{h}(s)), X^{\vv}_{h}(\eta_{h}(s-\tau)))|^p{\mbox d}s\\
&+C\vv^pt^{\frac{p}{2}-1}\mathbb{E}\int_0^t|g(X^{\vv}_{h}(\eta_{h}(s)), X^{\vv}_{h}(\eta_{h}(s-\tau)))|^p {\mbox d}s,
\end{split}
\end{equation}
where $Y_{h}^\vv(0)=\xi(0)-\theta f(\xi(0),\xi(-\tau))h$. Let $t=nh$ and $s=\bar{n}h$, where $n$ and $\bar{n}$ are nonnegative integers such that $\bar{n}h\le nh\le T$, then by assumption (H),
\begin{equation*}
\begin{split}
\E \left[\sup_{\bar{n}\le n}|Y_{h}^\vv(\bar{n}h)|^p\right]\le& C+C\sum\limits_{i=0}^{n-1}\E\left[\sup_{\bar{n}\le i}|X^{\vv}_{h}(\bar{n}h)|^p\right]h.
\end{split}
\end{equation*}
Since $\theta h<\frac{1}{6\beta}$, by $|x-y|^p\ge 2^{1-p}|x|^p-|y|^p$, we have
\begin{equation*}
\begin{split}
|Y_{h}^\vv(\bar{n}h)|^p\ge2^{1-p}|X_{h}^\vv(\bar{n}h)|^p-3^{p-1}\beta^p\theta^ph^p(1+|X_{h}^\vv(\bar{n}h)|^p+|X_{h}^\vv(\bar{n}h-mh)|^p),
\end{split}
\end{equation*}
which implies that
\begin{equation}\label{y0}
\begin{split}
\E \left[\sup_{\bar{n}\le n}|X_{h}^\vv(\bar{n}h)|^p\right]\le& C+C\E \left[\sup_{\bar{n}\le n}|Y_{h}^\vv(\bar{n}h)|^p\right]\\
\le& C+C\sum\limits_{i=0}^{n-1}\E\left[\sup_{\bar{n}\le i}|X^{\vv}_{h}(\bar{n}h)|^p\right]h.
\end{split}
\end{equation}
By the discrete Gronwall inequality,
\begin{equation}\label{disc}
\begin{split}
\E \left[\sup_{\bar{n}\le n}|X_{h}^\vv(\bar{n}h)|^p\right]\le C.
\end{split}
\end{equation}
Furthermore, by \eqref{july} and \eqref{disc},
\begin{equation*}
\begin{split}
\E \left[\sup_{0\le t\le T}|Y_{h}^\vv(t)|^p\right]\le& C+C\E\int_0^T(|X^{\vv}_{h}(\eta_{h}(s))|^p+|X^{\vv}_{h}(\eta_{h}(s-\tau))|^p){\mbox d}s\le C.
\end{split}
\end{equation*}
In the same way as \eqref{y0}, we derive
\begin{equation*}
\begin{split}
\E \left[\sup_{0\le t\le T}|X_{h}^\vv(t)|^p\right]\le& C+C\E \left[\sup_{0\le t\le T}|Y_{h}^\vv(t)|^p\right]\le C.
\end{split}
\end{equation*}
This completes the proof.\hfill $\Box$

\begin{lemma}\label{womiss}
Let assumption (H) hold. Then, for any $p\ge2,$ we have
\begin{equation*}
\sup_{0\le n\le M-1}\E [|X^\vv_{h}(t_{n+1})-X_{h}^\vv(t_n)|^p]\le Ch^p+C \vv^ph^{p/2}.
\end{equation*}
\end{lemma}
{\bf Proof.} Use the same notation $Y_h^\vv(t)$ as in Lemma \ref{0pmoment}. We derive from \eqref{continuous}, assumption (H), Lemma \ref{0pmoment} and the BDG inequality that for $p\ge 2$
\begin{equation*}
\begin{split}
\E|Y^\vv_{h}(t_{n+1})-Y_{h}^\vv(t_n)|^p\le&2^{p-1}h^{p-1}\int_{t_n}^{t_{n+1}}|f(X_h^{\vv}(\eta_h(s)), X_h^{\vv}(\eta_h(s-\tau)))|^p{\mbox d}s\\
&+\vv^ph^{\frac{p}{2}-1}\int_{t_n}^{t_{n+1}}|g(X_h^{\vv}(\eta_h(s)), X_h^{\vv}(\eta_h(s-\tau)))|^p{\mbox d}s\\
\le& Ch^p+C\vv^ph^{\frac{p}{2}}.
\end{split}
\end{equation*}
With the relationship between $X_h^\vv(t)$ and $Y_h^\vv(t)$, we obtain
\begin{equation*}
\begin{split}
X^\vv_{h}(t_{n+1})-X_{h}^\vv(t_n)=&Y^\vv_{h}(t_{n+1})-Y_{h}^\vv(t_n)+\theta f(X_h^\vv(t_{n+1}),X_h^\vv(t_{n+1}-\tau))h\\
&-\theta f(X_h^\vv(t_n),X_h^\vv(t_n-\tau))h.
\end{split}
\end{equation*}
By Lemma \ref{0pmoment}, it is easy to show that
\begin{equation*}
\begin{split}
\E|X^\vv_{h}(t_{n+1})-X_{h}^\vv(t_n)|^p\le C\E|Y^\vv_{h}(t_{n+1})-Y_{h}^\vv(t_n)|^p+Ch^p\le  Ch^p+C\vv^ph^{\frac{p}{2}}.
\end{split}
\end{equation*}
\hfill $\Box$

We now reveal the error between the numerical solution \eqref{continuous} and the exact solution \eqref{2.0}.

\begin{theorem}\label{th0}
Let assumption (H) hold, assume that $\Psi:\mathbb{R}^a\rightarrow\mathbb{R}$ has continuous second order derivative and there exists a constant $C$ such that
\begin{equation*}
\begin{split}
\left|\frac{\partial \Psi}{\partial x_i}\right|\le C~~and ~~\left|\frac{\partial^2 \Psi}{\partial x_i\partial x_j}\right|\le C
\end{split}
\end{equation*}
for any $i,j=1,2,\cdots,a$. Then, we have
\begin{equation*}
\begin{split}
\sup\limits_{0\le t\le T}\mathbb{E}|\Psi(X^\vv(t))-\Psi(X_h^\vv(t))|^2
=\mathcal{O}(h^2+\vv^2h).
\end{split}
\end{equation*}
\end{theorem}
{\bf Proof. } Set $I(t)=X_h^\vv(t)-X^{\vv}(t)-\theta f(X_h^\vv(t),X_h^\vv(t-\tau))h$, then
\begin{equation*}
\begin{split}
I(t)&=I(0)+\int_0^t[f(X_h^{\vv}(\eta_h(s)), X_h^{\vv}(\eta_h(s-\tau)))-f(X^{\vv}(s), X^{\vv}(s-\tau))]{\mbox d}s\\
&+\vv\int_0^t[g(X_h^{\vv}(\eta_h(s)), X_h^{\vv}(\eta_h(s-\tau)))-g(X^{\vv}(s), X^{\vv}(s-\tau))]{\mbox d}W(s),
\end{split}
\end{equation*}
where $I(0)=-\theta f(\xi(0),\xi(-\tau))h$. By the assumption (H) and Lemma \ref{womiss}, we see
\begin{equation}\label{new}
\begin{split}
\sup\limits_{0\le t\le T}\mathbb{E}|I(t)|^2&\le3|I(0)|^2+3T\mathbb{E}\int_0^T|f(X_h^{\vv}(\eta_h(s)), X_h^{\vv}(\eta_h(s-\tau)))-f(X^{\vv}(s), X^{\vv}(s-\tau))|^2{\mbox d}s\\
&+C\vv^2\mathbb{E}\int_0^T|g(X_h^{\vv}(\eta_h(s)), X_h^{\vv}(\eta_h(s-\tau)))-g(X^{\vv}(s), X^{\vv}(s-\tau))|^2{\mbox d}s\\
&\le Ch^2+ C\vv^2h+C\vv^2\int_0^T\sup\limits_{0\le t\le s}\mathbb{E}|X^{\vv}(t)-X_h^{\vv}(t)|^2{\mbox d}s+C\vv^2 h^2+C\vv^4h.
\end{split}
\end{equation}
Since $|x-y|^p\ge\frac{1}{2}|x|^p-|y|^p$, we get
\begin{equation*}
\begin{split}
|I(t)|^2\ge\frac{1}{2}|X^\vv(t)-X^\vv_h(t)|^2-|\theta f(X_h^\vv(t),X_h^\vv(t-\tau))h|^2
\end{split}
\end{equation*}
We then derive from Lemma \ref{0pmoment}, \eqref{new} and the Gronwall inequality that
\begin{equation*}
\begin{split}
\sup\limits_{0\le t\le T}\mathbb{E}|X^\vv(t)-X_h^\vv(t)|^2\le Ch^2+\vv^2h.
\end{split}
\end{equation*}
Since $\Psi$ has continuous bounded first order derivative, we immediately get
\begin{equation*}
\begin{split}
\sup\limits_{0\le t\le T}\mathbb{E}|\Psi(X^\vv(t))-\Psi(X_h^\vv(t))|^2\le C\sup\limits_{0\le t\le T}\mathbb{E}|X^\vv(t)-X_h^\vv(t)|^2.
\end{split}
\end{equation*}
The desired result then follows. \hfill $\Box$

\begin{coro}
Assume that the conditions of Theorem \ref{th0} hold. Let $M\ge 2, l\ge 1$, $h_l=T\cdot M^{-l}, h_{l-1}=T\cdot M^{-(l-1)}$.  Then
\begin{equation*}
\begin{split}
\sup_{0\le n<M^{l-1}}{\rm Var}(\Psi(X_{h_l}^\vv(t_n))-\Psi(X_{h_{l-1}}^\vv(t_n)))\le Ch_{l-1}^2+C\vv^2h_{l-1}.
\end{split}
\end{equation*}
\end{coro}
{\bf Proof.} For $0\le n\le M^{l-1}-1$, by Theorem \ref{th0},
\begin{equation*}
\begin{split}
&{\rm Var}(\Psi(X_{h_l}^\vv(t_n))-\Psi(X_{h_{l-1}}^\vv(t_n)))\le\mathbb{E}|\Psi(X_{h_l}^\vv(t_n))-\Psi(X_{h_{l-1}}^\vv(t_n))|^2\\
\le&2\mathbb{E}|\Psi(X_{h_l}^\vv(t_n))-\Psi(X^\vv(t_n))|^2+2\mathbb{E}|\Psi(X^\vv(t_n))-\Psi(X_{h_{l-1}}^\vv(t_n))|^2\\
\le&Ch_{l-1}^2+C\vv^2h_{l-1}.
\end{split}
\end{equation*}
\hfill $\Box$

\subsection{The Multilevel Monte Carlo theta EM Scheme}\label{sec2.2}
We now define the multilevel Monte Carlo theta EM scheme. Given any $T>0$, let $M\ge 2, l\ge 1$, $h_l=T\cdot M^{-l}, h_{l-1}=T\cdot M^{-(l-1)},$ assume there exists an $m_l$ such that $\tau=m_lh_l$. Let
\begin{equation}\label{c1}
\begin{split}
&X^{\vv}_{h_l}(t)-\theta f(X^{\vv}_{h_l}(t),X^{\vv}_{h_l}(t-\tau))h_l\\
=&\xi(0)-\theta f(\xi(0),\xi(-\tau))h_l+\int_0^tf(X^{\vv}_{h_l}(\eta_{h_{l}}(s)), X^{\vv}_{h_l}(\eta_{h_{l}}(s-\tau))){\mbox d}s\\
&+\vv\int_0^tg(X^{\vv}_{h_l}(\eta_{h_{l}}(s)), X^{\vv}_{h_l}(\eta_{h_{l}}(s-\tau))){\mbox d}W(s),
\end{split}
\end{equation}
and
\begin{equation}\label{c2}
\begin{split}
&X^{\vv}_{h_{l-1}}(t)-\theta f(X^{\vv}_{h_{l-1}}(t),X^{\vv}_{h_{l-1}}(t-\tau))h_{l-1}\\
=&\xi(0)-\theta f(\xi(0),\xi(-\tau))h_{l-1}+\int_0^tf(X^{\vv}_{h_{l-1}}(\eta_{h_{l-1}}(s)), X^{\vv}_{h_{l-1}}(\eta_{h_{l-1}}(s-\tau))){\mbox d}s\\
&+\vv\int_0^tg(X^{\vv}_{h_{l-1}}(\eta_{h_{l-1}}(s)), X^{\vv}_{h_{l-1}}(\eta_{h_{l-1}}(s-\tau))){\mbox d}W(s),
\end{split}
\end{equation}
where $\eta_{h_l}(s)=\left\lfloor s/h_l\right\rfloor h_l$. Here $\theta\in[0,1]$ is a parameter to control the implicitness. For $n\in \{0, 1, \ldots, M^{l-1}-1\}$ and $k\in \{0, \ldots, M\}$, let
$$
t_n=nh_{l-1} \mbox{ and } t_n^k=nh_{l-1}+kh_l.
$$
This means we divide the interval $[t_n, t_{n+1}]$ into $M$ equal parts, we have $t_n^0=t_n, t_n^{M}=t_{n+1}.$ We can rewrite \eqref{c1} and \eqref{c2} as the following discretization schemes.
For  $n\in \{0, 1, \ldots, M^{l-1}-1\}$ and $k\in \{0, \ldots, M-1\}$, let
\begin{equation}\label{c3}
\begin{split}
&X^{\vv}_{h_l}(t_n^{k+1})-\theta f(X^{\vv}_{h_{l}}(t_n^{k+1}),X^{\vv}_{h_{l}}(t_{n}^{k+1}-m_lh_l))h_l\\
=& X^{\vv}_{h_l}(t_n^{k})-\theta f(X^{\vv}_{h_{l}}(t_n^{k}),X^{\vv}_{h_{l}}(t_{n}^{k}-m_lh_l))h_l+f(X^{\vv}_{h_l}(t_n^k), X^{\vv}_{h_l}(t_{n}^{k}-m_lh_l))h_l\\
&+\vv \sqrt{h_l}g (X^{\vv}_{h_l}(t_n^k), X^{\vv}_{h_l}(t_{n}^{k}-m_lh_l))\xi_n^k,
\end{split}
\end{equation}
where the random vector $\xi_n^k\in\mathbb{R}^d$ has independent components, and each component is distributed as $N(0, 1).$ This implies
\begin{equation}\label{c4}
\begin{split}
&X^{\vv}_{h_l}(t_{n+1})-\theta f(X^{\vv}_{h_{l}}(t_{n+1}),X^{\vv}_{h_{l}}(t_{n+1}-m_lh_l))h_l\\
=&X^{\vv}_{h_l}(t_n)-\theta f(X^{\vv}_{h_{l}}(t_{n}),X^{\vv}_{h_{l}}(t_{n}-m_lh_l))h_l+\sum_{k=0}^{M-1}f(X^{\vv}_{h_l}(t_n^k), X^{\vv}_{h_l}(t_{n}^{k}-m_lh_l))h_l\\
&+\vv \sqrt{h_l}\sum_{k=0}^{M-1}g (X^{\vv}_{h_l}(t_n^k), X^{\vv}_{h_l}(t_{n}^{k}-m_lh_l))\xi_n^k.
\end{split}
\end{equation}
To simulate $X_{h_{l-1}}^\vv,$ we use
\begin{equation}\label{c5}
\begin{split}
&X^{\vv}_{h_{l-1}}(t_{n+1})-\theta f(X^{\vv}_{h_{l-1}}(t_{n+1}),X^{\vv}_{h_{l-1}}(t_{n+1}-m_lh_l))h_{l-1}\\
=&X^{\vv}_{h_{l-1}}(t_n)-\theta f(X^{\vv}_{h_{l-1}}(t_{n}),X^{\vv}_{h_{l-1}}(t_{n}-m_lh_l))h_{l-1}+f(X^{\vv}_{h_{l-1}}(t_n), X^{\vv}_{h_{l-1}}(t_{n}-m_lh_l))h_{l-1}\\
&+\vv \sqrt{h_{l}}g (X^{\vv}_{h_{l-1}}(t_n), X^{\vv}_{h_{l-1}}(t_{n}-m_lh_l))\sum_{k=0}^{M-1}\xi_n^k.
\end{split}
\end{equation}
For convenience, let
\begin{equation*}
Y^{\vv}_{h_l}(t):=X^{\vv}_{h_l}(t)-\theta f(X^{\vv}_{h_{l}}(t),X^{\vv}_{h_{l}}(t-m_lh_l))h_l,
\end{equation*}
 and
\begin{equation*}
 Y^{\vv}_{h_{l-1}}(t):=X^{\vv}_{h_{l-1}}(t)-\theta f(X^{\vv}_{h_{l-1}}(t),X^{\vv}_{h_{l-1}}(t-m_lh_l))h_{l-1}.
\end{equation*}
We now have the following estimates.

\begin{lemma}\label{pmoment}
Let assumption (H) hold. Then, for any $T>0$ and $p\ge 2,$ we have
$$
\E \left[\sup_{0\le t\le T}|X_{h_{l}}^\vv(t)|^p\right] \le C,
$$
and
$$
\E \left[\sup_{0\le t\le T}|X_{h_{l-1}}^\vv(t)|^p\right] \le C.
$$
\end{lemma}
{\bf Proof.} We omit the proof here since it is similar to that of Lemma \ref{0pmoment}. \hfill $\Box$

Let $Z_h$ be the deterministic solution to
\begin{equation*}
Z_h(t)-\theta f(Z_h(t), Z_h(t-\tau))h=\xi(0)-\theta f(\xi(0),\xi(-\tau))h+\int_0^tf(Z_h(\eta_h(s)), Z_h(\eta_h(s-\tau))){\mbox d}s,
\end{equation*}
which is the theta EM approximation to the ordinary differential delay equation obtained from \eqref{2.0} by taking $\vv=0$.

\begin{lemma}\label{sun}
Let assumption (H) hold. Then, for any $T>0$ and $p\ge 2,$ we have
$$
\E \left[\sup_{0\le t\le T}|Z_{h_{l}}(t)|^p\right] \le C,
$$
and
$$
\sup_{0\le n<M^{l-1}, 1\le k\le M}\E |Z_{h_{l}}(t_n^k)-Z_{h_l}(t_n)|^p \le CM^ph_l^p.
$$
\end{lemma}
{\bf Proof.} Following the proof of Lemma \ref{0pmoment}, under global Lipschitz condition (H), the first part is obvious. Denote by $\bar{Z}_{h_l}(t)=Z_{h_l}(t)-\theta f(Z_{h_l}(t), Z_{h_l}(t-\tau))h_l$. By the result of the first part,
\begin{equation*}
\begin{split}
\E|\bar{Z}_{h_l}(t_n^k)-\bar{Z}_{h_l}(t_n)|^p\le |kh_l|^{p-1}\E\int_{t_n}^{t_n^k}|f(Z_{h_l}(\eta_{h_l}(s)), Z_{h_l}(\eta_{h_l}(s-\tau)))|^p{\mbox d}s\le CM^ph_l^p.
\end{split}
\end{equation*}
On the other side, we see
\begin{equation*}
\begin{split}
\E|Z_{h_l}(t_n^k)-Z_{h_l}(t_n)|^p\le C\E|\bar{Z}_{h_l}(t_n^k)-\bar{Z}_{h_l}(t_n)|^p+CM^ph_l^p.
\end{split}
\end{equation*}
Thus, the desired assertion follows. \hfill $\Box$

\begin{lemma}\label{tech}
Let assumption (H) hold. Then, for any $T>0$ and $p\ge 2,$ we have
$$
\E \left[\sup_{0\le t\le T}|X_{h_{l}}^\vv(t)-Z_{h_l}(t)|^p\right] \le C \vv^p,
$$
and
$$
\E \left[\sup_{0\le t\le T}|X_{h_{l-1}}^\vv(t)-Z_{h_{l-1}}(t)|^p\right] \le C \vv^p.
$$
\end{lemma}
{\bf Proof.} Use the notation $\bar{Z}_{h_l}(t)$ defined in Lemma \ref{sun}. By the definition of $Y^{\vv}_{h_l}(t)$ and $\bar{Z}_{h_l}(t)$,
\begin{equation*}
\begin{split}
Y^{\vv}_{h_l}(t)-\bar{Z}_{h_l}(t)=&\int_0^t[f(X^{\vv}_{h_l}(\eta_{h_{l}}(s)), X^{\vv}_{h_l}(\eta_{h_{l}}(s-\tau)))-f(Z_{h_l}(\eta_{h_l}(s)), Z_{h_l}(\eta_{h_l}(s-\tau)))]{\mbox d}s\\
&+\vv\int_0^tg(X^{\vv}_{h_l}(\eta_{h_{l}}(s)), X^{\vv}_{h_l}(\eta_{h_{l}}(s-\tau))){\mbox d}W(s),
\end{split}
\end{equation*}
thus, by the BDG inequality, we get
\begin{equation}\label{guo}
\begin{split}
&\E\left[\sup\limits_{0\le s\le t}|Y^{\vv}_{h_l}(s)-\bar{Z}_{h_l}(s)|^p\right]\\
\le&C\int_0^t|f(X^{\vv}_{h_l}(\eta_{h_{l}}(s)), X^{\vv}_{h_l}(\eta_{h_{l}}(s-\tau)))-f(Z_{h_l}(\eta_{h_l}(s)), Z_{h_l}(\eta_{h_l}(s-\tau)))|^p{\mbox d}s\\
&+C\vv^p\int_0^t|g(X^{\vv}_{h_l}(\eta_{h_{l}}(s)), X^{\vv}_{h_l}(\eta_{h_{l}}(s-\tau)))|^p{\mbox d}s.
\end{split}
\end{equation}
Let $t=nh_l$ and $s=\bar{n}h_l$, where $n$ and $\bar{n}$ are nonnegative integers such that $\bar{n}h_l\le nh_l \le T$, then by assumption (H),
\begin{equation}\label{yuan}
\begin{split}
\E\left[\sup\limits_{\bar{n}\le n}|Y^{\vv}_{h_l}(\bar{n}h_l)-\bar{Z}_{h_l}(\bar{n}h_l)|^p\right]\le&C\sum\limits_{i=0}^{n-1}\E\left[\sup\limits_{\bar{n}\le i}|X^{\vv}_{h_l}(\bar{n}h_l)-Z_{h_l}(\bar{n}h_l)|^p\right]h_l\\
&+C\vv^p+C\vv^p\E\left[\sup\limits_{0\le s\le t}|X^{\vv}_{h_l}(s)|^p\right].
\end{split}
\end{equation}
By using $|x-y|^p\ge 2^{1-p}|x|^p-|y|^p$ and assumption (H) again, we see
\begin{equation*}
\begin{split}
\E\left[\sup\limits_{\bar{n}\le n}|X^{\vv}_{h_l}(\bar{n}h_l)-Z_{h_l}(\bar{n}h_l)|^p\right]\le C\E\left[\sup\limits_{\bar{n}\le n}|Y^{\vv}_{h_l}(\bar{n}h_l)-\bar{Z}_{h_l}(\bar{n}h_l)|^p\right],
\end{split}
\end{equation*}
then, Lemma \ref{pmoment} and \eqref{yuan} give that
\begin{equation*}
\begin{split}
&\E\left[\sup\limits_{\bar{n}\le n}|X^{\vv}_{h_l}(\bar{n}h_l)-Z_{h_l}(\bar{n}h_l)|^p\right]\\
\le&C\sum\limits_{i=0}^{n-1}\E\left[\sup\limits_{\bar{n}\le i}|X^{\vv}_{h_l}(\bar{n}h_l)-Z_{h_l}(\bar{n}h_l)|^p\right]h_l+C\vv^p.
\end{split}
\end{equation*}
The discrete Gronwall inequality leads to
\begin{equation*}
\begin{split}
&\E\left[\sup\limits_{\bar{n}\le n}|X^{\vv}_{h_l}(\bar{n}h_l)-Z_{h_l}(\bar{n}h_l)|^p\right]\le C\vv^p.
\end{split}
\end{equation*}
Furthermore, with assumption (H), we derive from \eqref{guo} and Lemma \ref{pmoment} that
\begin{equation*}
\begin{split}
&\E\left[\sup\limits_{0\le t\le T}|Y^{\vv}_{h_l}(t)-\bar{Z}_{h_l}(t)|^p\right]\le C\vv^p.
\end{split}
\end{equation*}
Then, the first part follows by using the relationship between $Y^{\vv}_{h_l}(t)$ and $X^{\vv}_{h_l}(t)$ together with Lemma \ref{pmoment}. By the same technique, the second part can be verified. \hfill $\Box$

\begin{lemma}\label{jia}
Let assumption (H) hold. Then, we have
\begin{equation*}
\sup_{0\le n<M^{l-1}, 1\le k\le M}|\E [X^\vv_{h_{l}}(t_n^k)-X_{h_l}^\vv(t_n)]|\le CMh_{l}.
\end{equation*}
\end{lemma}
{\bf Proof.} By \eqref{c3}, for $k\in\{1,2,\cdots,M\}$,
\begin{equation}\label{haha}
\begin{split}
Y^{\vv}_{h_l}(t_{n}^k)=&Y^{\vv}_{h_l}(t_{n}^0)+\sum_{q=0}^{k-1}f(X^{\vv}_{h_l}(t_n^q), X^{\vv}_{h_l}(t_{n}^{q}-m_lh_l))h_l\\
&+\vv \sqrt{h_l}\sum_{q=0}^{k-1}g (X^{\vv}_{h_l}(t_n^q), X^{\vv}_{h_l}(t_{n}^{q}-m_lh_l))\xi_n^q.
\end{split}
\end{equation}
Taking expectation on both sides, together with Lemma \ref{pmoment}, yields
\begin{equation}\label{one}
\begin{split}
|\E[Y^{\vv}_{h_l}(t_{n}^k)-Y^{\vv}_{h_l}(t_{n}^0)]|\le&\left|\E\left[\sum_{q=0}^{k-1}f(X^{\vv}_{h_l}(t_n^q), X^{\vv}_{h_l}(t_{n}^{q}-m_lh_l))h_l\right]\right|\\
&+\left|\E\left[\vv\sqrt{h_l}\sum_{q=0}^{k-1}g (X^{\vv}_{h_l}(t_n^q), X^{\vv}_{h_l}(t_{n}^{q}-m_lh_l))\xi_n^q\right]\right|\\
%\le&h_l\sum_{i=0}^{k-1}\E\left[\left|f(X^{\vv}_{h_l}(t_n^i), X^{\vv}_{h_l}(t_{n}^{i}-m_lh_l))\right|\right]\\
\le&Ch_l\sum_{q=0}^{k-1}\left(1+\E|X^{\vv}_{h_l}(t_n^q)|+\E|X^{\vv}_{h_l}(t_{n}^{q}-m_lh_l)|\right)\\
\le&CMh_{l}.
\end{split}
\end{equation}
Since we have
\begin{equation}\label{two}
\begin{split}
X^{\vv}_{h_l}(t_{n}^k)-X^{\vv}_{h_l}(t_{n})=&Y^{\vv}_{h_l}(t_{n}^k)-Y^{\vv}_{h_l}(t_{n}^0)+\theta f(X^{\vv}_{h_{l}}(t_{n}^k),X^{\vv}_{h_{l}}(t_{n}^k-m_lh_l))h_l\\
&-\theta f(X^{\vv}_{h_{l}}(t_{n}^0),X^{\vv}_{h_{l}}(t_{n}^0-m_lh_l))h_l.
\end{split}
\end{equation}
Combining \eqref{one} and \eqref{two}, it is easy to show the desired result by Lemma \ref{pmoment}.\hfill $\Box$

\begin{lemma}\label{miss}
Let assumption (H) hold. Then, for any $p>0,$ we have
\begin{equation*}
\sup_{0\le n<M^{l-1}, 1\le k\le M}\E [|X^\vv_{h_{l}}(t_n^k)-X_{h_l}^\vv(t_n)|^p]\le CM^p h_l^p+C \vv^pM^{p/2}h_l^{p/2}.
\end{equation*}
\end{lemma}
{\bf Proof.} We derive from \eqref{haha}, assumption (H), Lemma \ref{pmoment} and the discrete BDG inequality that for $p\ge 2$
\begin{equation*}
\begin{split}
\E|Y^{\vv}_{h_l}(t_{n}^k)-Y^{\vv}_{h_l}(t_{n}^0)|^p\le&2^{p-1}M^{p-1}h_l^p\sum_{q=0}^{k-1}\E|f(X^{\vv}_{h_l}(t_n^q), X^{\vv}_{h_l}(t_{n}^{q}-m_lh_l))|^p\\
&+2^{p-1}\vv^ph_l^{p/2}\E\left|\sum_{q=0}^{k-1}g (X^{\vv}_{h_l}(t_n^q), X^{\vv}_{h_l}(t_{n}^{q}-m_lh_l))\xi_n^q\right|^p\\
\le&CM^ph_l^p+C\vv^ph_l^{p/2}\E\left|\sum\limits_{q=0}^{k-1}\sum\limits_{i=1}^{d}|g^i(X^{\vv}_{h_l}(t_n^q), X^{\vv}_{h_l}(t_{n}^{q}-m_lh_l))|^2\right|^{\frac{p}{2}}\\
\le&CM^ph_l^p+C\vv^p M^{p/2}h_l^{p/2}\left(1+\sup_{0\le q\le M-1}\E|X^{\vv}_{h_l}(t_n^q)|^p\right)\\
\le&CM^ph_l^p+C\vv^p M^{p/2}h_l^{p/2},
\end{split}
\end{equation*}
where $g^i$ is the $i$-th column of $g$ and in the last step we have used Lemma \ref{pmoment}. By \eqref{two} and Lemma \ref{pmoment} again, we conclude that the desired result follows for $p\ge2$. Finally, one can use the Young inequality to get the results for $p\in(0,2)$.\hfill $\Box$

Taylor expansion of the drift coefficient.
\begin{lemma}\label{taylor}
Let $\nabla$ and $\nabla^2$ be the first and second order derivatives respectively. Then
\begin{equation*}
\begin{split}
&f(X^\vv_{h_l}(t_n^k), X^\vv_{h_l}(t_{n}^k-m_lh_l))-f(X^\vv_{h_l}(t_n), X^\vv_{h_l}(t_{n}-m_lh_l))=A_{n}^{k}+B_{n}^{k}+C_{n}^{k},
\end{split}
\end{equation*}
where $A_{n}^{k}, B_{n}^{k}, C_{n}^{k}$ are defined as in the proof.
\end{lemma}

\noindent
{\bf Proof.}
Let $f_i(x)$ be the $i$th component of $f(x)$. By the Taylor expansion, for $i=1,2,\cdots, a$,
\begin{equation*}
\begin{split}
&f_i(X^\vv_{h_l}(t_n^k), X^\vv_{h_l}(t_{n}^k-m_lh_l))-f_i(X^\vv_{h_l}(t_n), X^\vv_{h_l}(t_{n}-m_lh_l))\\
=& \int_0^1\nabla f_i\big((X^\vv_{h_l}(t_n), X^\vv_{h_l}(t_{n}-m_lh_l)) +s[(X^\vv_{h_l}(t_n^k), X^\vv_{h_l}(t_{n}^k-m_lh_l))\\
&-(X^\vv_{h_l}(t_n), X^\vv_{h_l}(t_{n}-m_lh_l))]\big)ds\times \begin{pmatrix}
X^\vv_{h_l}(t_n^k)-X^\vv_{h_l}(t_n)\\
X^\vv_{h_l}(t_{n}^k-m_lh_l)-X^\vv_{h_l}(t_{n}-m_lh_l)
\end{pmatrix}\\
=&
\int_0^1\nabla f_i\big((X^\vv_{h_l}(t_n), X^\vv_{h_l}(t_{n}-m_lh_l)) +s[(X^\vv_{h_l}(t_n^k), X^\vv_{h_l}(t_{n}^k-m_lh_l))\\
&-(X^\vv_{h_l}(t_n), X^\vv_{h_l}(t_{n}-m_lh_l))]\big)ds\times \begin{pmatrix}\sigma_n^{11}+\sigma_n^{12}+\sigma_n^{13}\\
\sigma_n^{21}+\sigma_n^{22}+\sigma_n^{23}
\end{pmatrix},
\end{split}
\end{equation*}
where
\begin{equation*}
\begin{split}
\sigma_n^{11}=&\sum_{q=0}^{k-1}f(X^{\vv}_{h_l}(t_n^q), X^{\vv}_{h_l}(t_{n}^q-m_lh_l))h_l,\\
\sigma_n^{12}=&\vv \sqrt{h_l}\sum_{q=0}^{k-1}g (X^{\vv}_{h_l}(t_n^q), X^{\vv}_{h_l}(t_{n}^q-m_lh_l)) \xi_n^q,\\
\sigma_n^{13}=&\theta f(X^{\vv}_{h_{l}}(t_{n}^k),X^{\vv}_{h_{l}}(t_{n}^k-m_lh_l))h_l-\theta f(X^{\vv}_{h_{l}}(t_{n}),X^{\vv}_{h_{l}}(t_{n}-m_lh_l))h_l, \\
\sigma_n^{21}=&\sum_{q=0}^{k-1}f(X^{\vv}_{h_l}(t_{n}^q-m_lh_l), X^{\vv}_{h_l}(t_{n}^q-2m_lh_l))h_l,\\
\sigma_n^{22}=&\vv \sqrt{h_l}\sum_{q=0}^{k-1}g (X^{\vv}_{h_l}(t_{n}^q-m_lh_l), X^{\vv}_{h_l}(t_{n}^q-2m_lh_l)) \xi_n^q,\\
\sigma_n^{23}=&\theta f(X^{\vv}_{h_{l}}(t_{n}^k-m_lh_l),X^{\vv}_{h_{l}}(t_{n}^k-2m_lh_l))h_l-\theta f(X^{\vv}_{h_{l}}(t_{n}-m_lh_l),X^{\vv}_{h_{l}}(t_{n}-2m_lh_l))h_l.
\end{split}
\end{equation*}
Again by the Taylor expansion we derive
\begin{equation*}
\begin{split}
 &\int_0^1\nabla f_i\big((X^\vv_{h_l}(t_n), X^\vv_{h_l}(t_{n}-m_lh_l)) +s[(X^\vv_{h_l}(t_n^k), X^\vv_{h_l}(t_{n}^k-m_lh_l))\\
&-(X^\vv_{h_l}(t_n), X^\vv_{h_l}(t_{n}-m_lh_l))]\big)ds\times \begin{pmatrix}
\sigma_n^{12}\\
\sigma_n^{22}
\end{pmatrix}\\
=&\nabla f_i(X^\vv_{h_l}(t_n), X^\vv_{h_l}(t_{n}-m_lh_l))\times \begin{pmatrix}
\sigma_n^{12}\\
\sigma_n^{22}
\end{pmatrix}+[X_{h_l}^\vv(t_n^k)-X_{h_l}^\vv(t_n), X_{h_l}^\vv(t_n^k-m_lh_l)-X_{h_l}^\vv(t_n-m_lh_l)]\\
&\cdot\int_0^1\int_0^s\nabla^2 f_i\big((X^\vv_{h_l}(t_n), X^\vv_{h_l}(t_{n}-m_lh_l))+u[(X^\vv_{h_l}(t_n^k), X^\vv_{h_l}(t_{n}^k-m_lh_l)) \\
&-(X^\vv_{h_l}(t_n), X^\vv_{h_l}(t_{n}-m_lh_l))]\big)duds\times \begin{pmatrix}
\sigma_n^{12}\\
\sigma_n^{22}
\end{pmatrix}.
\end{split}
\end{equation*}
These implies
\begin{equation*}
\begin{split}
&f_i(X^\vv_{h_l}(t_n^k), X^\vv_{h_l}(t_{n}^k-m_lh_l))-f_i(X^\vv_{h_l}(t_n), X^\vv_{h_l}(t_{n}-m_lh_l))\\
=& \int_0^1\nabla f_i\big((X^\vv_{h_l}(t_n), X^\vv_{h_l}(t_{n}-m_lh_l)) +s[(X^\vv_{h_l}(t_n^k), X^\vv_{h_l}(t_{n}^k-m_lh_l))\\
&-(X^\vv_{h_l}(t_n), X^\vv_{h_l}(t_{n}-m_lh_l))]\big)ds\times\begin{pmatrix}\sigma_n^{11}+\sigma_n^{13}\\
\sigma_n^{21}+\sigma_n^{23}
\end{pmatrix}
+\nabla f_i(X^\vv_{h_l}(t_n), X^\vv_{h_l}(t_{n}-m_lh_l))\times \begin{pmatrix}
\sigma_n^{12}\\
\sigma_n^{22}
\end{pmatrix}\\
&+[X_{h_l}^\vv(t_n^k)-X_{h_l}^\vv(t_n), X_{h_l}^\vv(t_n^k-m_lh_l)-X_{h_l}^\vv(t_n-m_lh_l)]\\
&\cdot\int_0^1\int_0^s\nabla^2 f_i\big((X^\vv_{h_l}(t_n), X^\vv_{h_l}(t_{n}-m_lh_l)) +u[(X^\vv_{h_l}(t_n^k), X^\vv_{h_l}(t_{n}^k-m_lh_l))\\
&-(X^\vv_{h_l}(t_n), X^\vv_{h_l}(t_{n}-m_lh_l))]\big)duds\times \begin{pmatrix}
\sigma_n^{12}\\
\sigma_n^{22}
\end{pmatrix}\\
=:&A_n^{ik}+B_n^{ik}+C_{n}^{ik}.
\end{split}
\end{equation*}
Denote by $A_n^k=(A_n^{1k}, A_n^{2k},\cdots, A_n^{ak})^T$, $B_n^k=(B_n^{1k}, B_n^{2k},\cdots, B_n^{ak})^T$ and $C_n^k=(C_{n}^{1k}, C_{n}^{2k},\cdots, C_n^{ak})^T$, then we can rewrite $f(X^\vv_{h_l}(t_n^k), X^\vv_{h_l}(t_{n}^k-m_lh_l))-f(X^\vv_{h_l}(t_n), X^\vv_{h_l}(t_{n}-m_lh_l))$ as $A_{n}^{k}+B_{n}^{k}+C_{n}^{k}$. This completes the proof.\hfill $\Box$

\begin{theorem}\label{election}
Let assumption (H) hold. Then we have
\begin{equation*}
\sup_{0\le n<M^{l-1}}\E [|X^\vv_{h_{l}}(t_n)-X_{h_{l-1}}^\vv(t_n)|^2]\le CM^2h_{l}^2+C\vv^4Mh_{l}.
\end{equation*}
\end{theorem}
{\bf Proof.} For any $n\le M^{l-1}-1$ , by \eqref{c4} and \eqref{c5}, we get
\begin{equation*}
\begin{split}
&Y^\vv_{h_{l}}(t_{n+1})-Y_{h_{l-1}}^\vv(t_{n+1})=Y^\vv_{h_{l}}(t_{n})-Y_{h_{l-1}}^\vv(t_{n})\\
&+h_l\sum_{k=0}^{M-1}[f(X^{\vv}_{h_l}(t_n^k), X^{\vv}_{h_l}(t_{n}^{k}-m_lh_l))-f(X^{\vv}_{h_{l-1}}(t_n), X^{\vv}_{h_{l-1}}(t_{n}-m_lh_l))]\\
&+\vv \sqrt{h_l}\sum_{k=0}^{M-1}[g(X^{\vv}_{h_l}(t_n^k), X^{\vv}_{h_l}(t_{n}^{k}-m_lh_l))-g(X^{\vv}_{h_{l-1}}(t_n), X^{\vv}_{h_{l-1}}(t_{n}-m_lh_l))]\xi_n^k\\
=& Y^\vv_{h_{l}}(t_{n})-Y_{h_{l-1}}^\vv(t_{n})+h_l\sum_{k=0}^{M-1}[f(X^{\vv}_{h_l}(t_n^k), X^{\vv}_{h_l}(t_{n}^{k}-m_lh_l))-f(X^{\vv}_{h_{l}}(t_n), X^{\vv}_{h_{l}}(t_{n}-m_lh_l))]\\
&+h_l\sum_{k=0}^{M-1}[f(X^{\vv}_{h_l}(t_n), X^{\vv}_{h_{l}}(t_{n}-m_lh_l))-f(X^{\vv}_{h_{l-1}}(t_n), X^{\vv}_{h_{l-1}}(t_{n}-m_lh_l))]\\
&+\vv \sqrt{h_l}\sum_{k=0}^{M-1}[g(X^{\vv}_{h_l}(t_n^k), X^{\vv}_{h_l}(t_{n}^{k}-m_lh_l))-g(X^{\vv}_{h_{l}}(t_n), X^{\vv}_{h_{l}}(t_{n}-m_lh_l))]\xi_n^k\\
&+\vv \sqrt{h_l}\sum_{k=0}^{M-1}[g(X^{\vv}_{h_l}(t_n), X^{\vv}_{h_l}(t_{n}-m_lh_l))-g(X^{\vv}_{h_{l-1}}(t_n), X^{\vv}_{h_{l-1}}(t_{n}-m_lh_l))]\xi_n^k.
\end{split}
\end{equation*}
By the elementary inequality $\left(\sum\limits_{i=1}^nx_i\right)^2\le x_1^2+(n-1)\sum\limits_{i=2}^n|x_i|^2+2\sum\limits_{i=2}^n\langle x_1, x_i\rangle$, we compute
\begin{equation*}
\begin{split}
&|Y^\vv_{h_{l}}(t_{n+1})-Y_{h_{l-1}}^\vv(t_{n+1})|^2\le |Y^\vv_{h_{l}}(t_{n})-Y_{h_{l-1}}^\vv(t_{n})|^2\\
&+4M h_l^2\sum_{k=0}^{M-1}\left|f(X^{\vv}_{h_l}(t_n^k), X^{\vv}_{h_l}(t_{n}^{k}-m_lh_l))-f(X^{\vv}_{h_{l}}(t_n), X^{\vv}_{h_{l}}(t_{n}-m_lh_l))\right|^2\\
&+4M h_l^2\sum_{k=0}^{M-1}\left|f(X^{\vv}_{h_l}(t_n), X^{\vv}_{h_{l}}(t_{n}-m_lh_l))-f(X^{\vv}_{h_{l-1}}(t_n), X^{\vv}_{h_{l-1}}(t_{n}-m_lh_l))\right|^2\\
&+4\vv^2h_l\left|\sum_{k=0}^{M-1}[g(X^{\vv}_{h_l}(t_n^k), X^{\vv}_{h_l}(t_{n}^{k}-m_lh_l))-g(X^{\vv}_{h_{l}}(t_n), X^{\vv}_{h_{l}}(t_{n}-m_lh_l))]\xi_n^k\right|^2\\
&+4\vv^2h_l\left|\sum_{k=0}^{M-1}[g(X^{\vv}_{h_l}(t_n), X^{\vv}_{h_l}(t_{n}-m_lh_l))-g(X^{\vv}_{h_{l-1}}(t_n), X^{\vv}_{h_{l-1}}(t_{n}-m_lh_l))]\xi_n^k\right|^2\\
&+2h_l\sum_{k=0}^{M-1}\langle Y^\vv_{h_{l}}(t_{n})-Y_{h_{l-1}}^\vv(t_{n}), f(X^{\vv}_{h_l}(t_n^k), X^{\vv}_{h_l}(t_{n}^{k}-m_lh_l))-f(X^{\vv}_{h_{l}}(t_n), X^{\vv}_{h_{l}}(t_{n}-m_lh_l))\rangle\\
&+2h_l\sum_{k=0}^{M-1}\langle Y^\vv_{h_{l}}(t_{n})-Y_{h_{l-1}}^\vv(t_{n}), f(X^{\vv}_{h_l}(t_n), X^{\vv}_{h_{l}}(t_{n}-m_lh_l))-f(X^{\vv}_{h_{l-1}}(t_n), X^{\vv}_{h_{l-1}}(t_{n}-m_lh_l))\rangle\\
&+2\vv \sqrt{h_l}\sum_{k=0}^{M-1}\langle Y^\vv_{h_{l}}(t_{n})-Y_{h_{l-1}}^\vv(t_{n}), [g(X^{\vv}_{h_l}(t_n^k), X^{\vv}_{h_l}(t_{n}^{k}-m_lh_l))-g(X^{\vv}_{h_{l}}(t_n), X^{\vv}_{h_{l}}(t_{n}-m_lh_l))]\xi_n^k\rangle\\
&+2\vv \sqrt{h_l}\sum_{k=0}^{M-1}\langle Y^\vv_{h_{l}}(t_{n})-Y_{h_{l-1}}^\vv(t_{n}), [g(X^{\vv}_{h_l}(t_n), X^{\vv}_{h_l}(t_{n}-m_lh_l))-g(X^{\vv}_{h_{l-1}}(t_n), X^{\vv}_{h_{l-1}}(t_{n}-m_lh_l))]\xi_n^k\rangle.
\end{split}
\end{equation*}
Taking expectation, then summing both sides, using assumption (H) and the Young inequality, we obtain that for $\Lambda\le M^{l-1}-1$
\begin{equation*}
\begin{split}
&\sup\limits_{0\le n\le \Lambda+1}\mathbb{E}|Y^\vv_{h_{l}}(t_{n})-Y_{h_{l-1}}^\vv(t_{n})|^2\\
\le&8\alpha^2M h_l^2\sum_{j=0}^{\Lambda}\sum_{k=0}^{M-1}\mathbb{E}\left(|X^{\vv}_{h_l}(t_j^k)-X^{\vv}_{h_{l}}(t_j)|^2+|X^{\vv}_{h_l}(t_{j}^{k}-m_lh_l)-X^{\vv}_{h_{l}}(t_{j}-m_lh_l)|^2\right)\\
&+8\alpha^2M h_l^2\sum_{j=0}^{\Lambda}\sum_{k=0}^{M-1}\mathbb{E}\left(|X^{\vv}_{h_l}(t_j)-X^{\vv}_{h_{l-1}}(t_j)|^2+|X^{\vv}_{h_{l}}(t_{j}-m_lh_l)- X^{\vv}_{h_{l-1}}(t_{j}-m_lh_l)|^2\right)\\
&+8\alpha^2\vv^2h_l\sum_{j=0}^{\Lambda}\sum_{k=0}^{M-1}\mathbb{E}\left(|X^{\vv}_{h_l}(t_j^k)-X^{\vv}_{h_{l}}(t_j)|^2+|X^{\vv}_{h_l}(t_{j}^{k}-m_lh_l)- X^{\vv}_{h_{l}}(t_{j}-m_lh_l)|^2\right)\\
&+8\alpha^2\vv^2h_l\sum_{j=0}^{\Lambda}\sum_{k=0}^{M-1}\mathbb{E}\left(|X^{\vv}_{h_l}(t_j)-X^{\vv}_{h_{l-1}}(t_j)|^2+|X^{\vv}_{h_l}(t_{j}-m_lh_l)- X^{\vv}_{h_{l-1}}(t_{j}-m_lh_l)|^2\right)\\
&+2h_l\sum_{j=0}^{\Lambda}\sum_{k=0}^{M-1}\mathbb{E}\langle Y^\vv_{h_{l}}(t_{j})-Y_{h_{l-1}}^\vv(t_{j}), f(X^{\vv}_{h_l}(t_j^k), X^{\vv}_{h_l}(t_{j}^{k}-m_lh_l))-f(X^{\vv}_{h_{l}}(t_j), X^{\vv}_{h_{l}}(t_{j}-m_lh_l))\rangle\\
&+2h_l\sum_{j=0}^{\Lambda}\sum_{k=0}^{M-1}\mathbb{E}\langle Y^\vv_{h_{l}}(t_{j})-Y_{h_{l-1}}^\vv(t_{j}), f(X^{\vv}_{h_l}(t_j), X^{\vv}_{h_{l}}(t_{j}-m_lh_l))-f(X^{\vv}_{h_{l-1}}(t_j), X^{\vv}_{h_{l-1}}(t_{j}-m_lh_l))\rangle\\
\end{split}
\end{equation*}
By Lemma \ref{miss}, we immediately get
\begin{equation}\label{bankho}
\begin{split}
&\sup\limits_{0\le n\le \Lambda+1}\mathbb{E}|Y^\vv_{h_{l}}(t_{n})-Y_{h_{l-1}}^\vv(t_{n})|^2\\
\le&CMh_l(M^2h_l^2+\vv^2Mh_l)+C\vv^2(M^2h_l^2+\vv^2Mh_l)\\
&+C(M^2 h_l^2+\vv^2M h_l)\sum_{j=0}^{\Lambda}\mathbb{E}|X^\vv_{h_{l}}(t_{j})-X_{h_{l-1}}^\vv(t_{j})|^2+\frac{1}{4}\sup\limits_{0\le n\le \Lambda+1}\mathbb{E}|Y^\vv_{h_{l}}(t_{n})-Y_{h_{l-1}}^\vv(t_{n})|^2\\
&+C(M^2 h_l^2+\vv^2M h_l)\sum_{j=0}^{\Lambda}\mathbb{E}|X^\vv_{h_{l}}(t_{j}-m_lh_l)-X_{h_{l-1}}^\vv(t_{j}-m_lh_l)|^2\\
&+2h_l\sum_{j=0}^{\Lambda}\sum_{k=0}^{M-1}\mathbb{E}\langle Y^\vv_{h_{l}}(t_{j})-Y_{h_{l-1}}^\vv(t_{j}), f(X^{\vv}_{h_l}(t_j^k), X^{\vv}_{h_l}(t_{j}^{k}-m_lh_l))-f(X^{\vv}_{h_{l}}(t_j), X^{\vv}_{h_{l}}(t_{j}-m_lh_l))\rangle\\
&+CM h_l\sum_{j=0}^{\Lambda}\mathbb{E}|X^\vv_{h_{l}}(t_{j})-X_{h_{l-1}}^\vv(t_{j})|^2+CM h_l\sum_{j=0}^{\Lambda}\mathbb{E}|X^\vv_{h_{l}}(t_{j}-m_lh_l)-X_{h_{l-1}}^\vv(t_{j}-m_lh_l)|^2.
\end{split}
\end{equation}
Applying the Young inequality and Lemma \ref{taylor}, we see
\begin{equation}\label{reject}
\begin{split}
&2h_l\sum_{j=0}^{\Lambda}\sum_{k=0}^{M-1}\mathbb{E}\langle Y^\vv_{h_{l}}(t_{j})-Y_{h_{l-1}}^\vv(t_{j}), f(X^{\vv}_{h_l}(t_j^k), X^{\vv}_{h_l}(t_{j}^{k}-m_lh_l))-f(X^{\vv}_{h_{l}}(t_j), X^{\vv}_{h_{l}}(t_{j}-m_lh_l))\rangle\\
=&2h_l\sum_{j=0}^{\Lambda}\sum_{k=0}^{M-1}\mathbb{E}\langle Y^\vv_{h_{l}}(t_{j})-Y_{h_{l-1}}^\vv(t_{j}), A_j^k\rangle+2h_l\sum_{j=0}^{\Lambda}\sum_{k=0}^{M-1}\mathbb{E}\langle Y^\vv_{h_{l}}(t_{j})-Y_{h_{l-1}}^\vv(t_{j}), B_j^k\rangle\\
&+2h_l\sum_{j=0}^{\Lambda}\sum_{k=0}^{M-1}\mathbb{E}\langle Y^\vv_{h_{l}}(t_{j})-Y_{h_{l-1}}^\vv(t_{j}), C_j^k\rangle\\
\le&\frac{1}{4}\sup\limits_{0\le n\le \Lambda+1}\mathbb{E}|Y^\vv_{h_{l}}(t_{n})-Y_{h_{l-1}}^\vv(t_{n})|^2+Ch_l\sum_{j=0}^{\Lambda}\sum_{k=0}^{M-1}\mathbb{E}|A_j^k|^2+Ch_l\sum_{j=0}^{\Lambda}\sum_{k=0}^{M-1}\mathbb{E}|C_j^k|^2\\
\le&\frac{1}{4}\sup\limits_{0\le n\le \Lambda+1}\mathbb{E}|Y^\vv_{h_{l}}(t_{n})-Y_{h_{l-1}}^\vv(t_{n})|^2+CM^2h_l^2+C\vv^2M^3h_l^3+C\vv^4M^2h_l^2,
\end{split}
\end{equation}
since by Lemma \ref{pmoment} and Lemma \ref{miss}, it is easy to see $\mathbb{E}|A_j^k|^2\le CM^2h_l^2+C\vv^2M^3h_l^3$, moreover, by the H\"{o}lder inequality and Lemmas \ref{pmoment}, \ref{miss}, we have
\begin{equation*}
\begin{split}
& \mathbb{E}|C_j^k|^2\le C\vv^2h_l\mathbb{E}\bigg(\left[|X_{h_l}^\vv(t_j^k)-X_{h_l}^\vv(t_j)|^2+|X_{h_l}^\vv(t_j^k-m_lh_l)-X_{h_l}^\vv(t_j-m_lh_l)|^2\right]\\
&\cdot\bigg[\sum_{q=0}^{k-1}|g (X^{\vv}_{h_l}(t_j^q), X^{\vv}_{h_l}(t_{j}^q-m_lh_l)) \xi_j^q|^2
+\sum_{q=0}^{k-1}|g (X^{\vv}_{h_l}(t_{j}^q-m_lh_l), X^{\vv}_{h_l}(t_{j}^q-2m_lh_l)) \xi_j^q|^2\bigg]\bigg)\\
\le& C\vv^2h_l\left(\mathbb{E}|X_{h_l}^\vv(t_j^k)-X_{h_l}^\vv(t_j)|^4+\mathbb{E}|X_{h_l}^\vv(t_j^k-m_lh_l)-X_{h_l}^\vv(t_j-m_lh_l)|^4\right)^{1/2}\\
&\cdot\left(\sum_{q=0}^{k-1}\mathbb{E}|g (X^{\vv}_{h_l}(t_j^q), X^{\vv}_{h_l}(t_{j}^q-m_lh_l)) \xi_j^q|^4
+\sum_{q=0}^{k-1}\mathbb{E}|g (X^{\vv}_{h_l}(t_{j}^q-m_lh_l), X^{\vv}_{h_l}(t_{j}^q-2m_lh_l)) \xi_j^q|^4\right)^{1/2}\\
\le&C\vv^2M^3h_l^3+C\vv^4M^2h_l^2.
\end{split}
\end{equation*}
By the definition of $Y^{\vv}_{h_l}(t_{n})$ and $Y^{\vv}_{h_{l-1}}(t_{n})$, we have
\begin{equation*}
\begin{split}
&Y^{\vv}_{h_l}(t_{n})-Y^{\vv}_{h_{l-1}}(t_{n})=X^{\vv}_{h_l}(t_{n})-X^{\vv}_{h_{l-1}}(t_{n})\\
&-\theta f(X^{\vv}_{h_{l}}(t_{n}),X^{\vv}_{h_{l}}(t_{n}-m_lh_l))h_l+\theta f(X^{\vv}_{h_{l-1}}(t_{n}),X^{\vv}_{h_{l-1}}(t_{n}-m_lh_l))h_{l-1}.
\end{split}
\end{equation*}
Taking advantage of the elementary equality $2(|a|^2+|b|^2)\ge |a-b|^2\ge|a|^2-|b|^2$, we get
\begin{equation*}
\begin{split}
&|Y^{\vv}_{h_l}(t_{n})-Y^{\vv}_{h_{l-1}}(t_{n})|^2\ge|X^{\vv}_{h_l}(t_{n})-X^{\vv}_{h_{l-1}}(t_{n})|^2\\
&-|\theta f(X^{\vv}_{h_{l}}(t_{n}),X^{\vv}_{h_{l}}(t_{n}-m_lh_l))h_{l}-\theta f(X^{\vv}_{h_{l-1}}(t_{n}),X^{\vv}_{h_{l-1}}(t_{n}-m_lh_l))h_{l-1}|^2\\
\ge&|X^{\vv}_{h_l}(t_{n})-X^{\vv}_{h_{l-1}}(t_{n})|^2-2\theta^2h_l^2|f(X^{\vv}_{h_{l}}(t_{n}),X^{\vv}_{h_{l}}(t_{n}-m_lh_l))|^2\\
&-2\theta^2h_{l-1}^2|f(X^{\vv}_{h_{l-1}}(t_{n}),X^{\vv}_{h_{l-1}}(t_{n}-m_lh_l))|^2.
\end{split}
\end{equation*}
This, together with Lemma \ref{pmoment} imply
\begin{equation}\label{aid}
\begin{split}
&\sup\limits_{0\le n\le \Lambda+1}\mathbb{E}|X^\vv_{h_{l}}(t_{n})-X_{h_{l-1}}^\vv(t_{n})|^2\le\sup\limits_{0\le n\le \Lambda+1}\mathbb{E}|Y^\vv_{h_{l}}(t_{n})-Y_{h_{l-1}}^\vv(t_{n})|^2+Ch_{l-1}^2.
\end{split}
\end{equation}
Combining \eqref{bankho}-\eqref{aid} yields
\begin{equation*}\label{bank}
\begin{split}
&\sup\limits_{0\le n\le \Lambda+1}\mathbb{E}|X^\vv_{h_{l}}(t_{n})-X_{h_{l-1}}^\vv(t_{n})|^2\\
\le& C(M^3h_l^3+\vv^2M^2h_l^2+\vv^4Mh_l+M^2h_l^2+\vv^2M^3h_l^3+\vv^4M^2h_l^2)\\
&+C(M^2 h_l^2+\vv^2Mh_l+Mh_l)\sum\limits_{j=0}^\Lambda\sup\limits_{0\le n\le j}\mathbb{E}|X^\vv_{h_{l}}(t_{n})-X_{h_{l-1}}^\vv(t_{n})|^2.
\end{split}
\end{equation*}
By the discrete Gronwall inequality, the desired result can be obtained since the dominant term above is of order $M^2h_{l}^2$ and $\vv^4Mh_{l}$.\hfill $\Box$

The following two lemmas are from \cite{ahs15}.
\begin{lemma}\label{tech1}
Suppose $X_1(t)$ and $X_2(t)$ are stochastic processes on $\mathbb{R}^a$ and that $x_1(t)$ and $x_2(t)$ are deterministic processes on $\mathbb{R}^a$. Further, suppose that
\begin{equation*}
\begin{split}
\sup_{t\le T}\mathbb{E}|X_1(t)-x_1(t)|^2\le C_1\vv^2, ~~\sup_{s\le T}\mathbb{E}|X_2(t)-x_2(t)|^2\le C_2\vv^2,
\end{split}
\end{equation*}
for some $C_1, C_2$ and any $\vv\in(0,1)$. Assume that $\Phi:\mathbb{R}^a\rightarrow\mathbb{R}$ is Lipschitz with Lipschitz constant $C_L$. Then
\begin{equation*}
\begin{split}
\sup_{t\le T}{\rm Var}\left(\int_0^1\Phi(X_2(t)+s(X_1(t)-X_2(t))){\rm d} s\right)\le C_L^2C_1C_2\vv^2.
\end{split}
\end{equation*}

\end{lemma}

\begin{lemma}\label{tech2}
Suppose that $A^{\vv h}$ and $B^{\vv h}$ are families of random variables determined by scaling parameters $\vv$ and $h$. Further, suppose that there are positive constants $C_1, C_2, C_3$ such that for any $\vv\in(0,1)$, the following three conditions hold: \\
{\rm (i)} ${\rm Var}(A^{\vv h})\le C_1\vv^2$ uniformly in $h$.\\
{\rm (ii)} $|A^{\vv h}|\le C_2$ uniformly in $h$.\\
{\rm (iii)} $|\mathbb{E}B^{\vv h}|\le C_3h$.\\
Then
\begin{equation*}
\begin{split}
{\rm Var}(A^{\vv h}B^{\vv h})\le 3C_3^2C_1\vv^2 h^2+15C_2^2{\rm Var}(B^{\vv h}).
\end{split}
\end{equation*}

\end{lemma}

\begin{theorem}\label{th1}
Let assumption (H) hold, assume that $\Psi:\mathbb{R}^a\rightarrow\mathbb{R}$ has continuous second order derivative and there exists a constant $C$ such that
\begin{equation*}
\begin{split}
\left|\frac{\partial \Psi}{\partial x_i}\right|\le C~~and ~~\left|\frac{\partial^2 \Psi}{\partial x_i\partial x_j}\right|\le C
\end{split}
\end{equation*}
for any $i,j=1,2,\cdots,a$. Then, we have
\begin{equation*}
\begin{split}
\sup_{0\le n<M^{l-1}}{\rm Var}(\Psi(X_{h_l}^\vv(t_n))-\Psi(X_{h_{l-1}}^\vv(t_n)))\le Ch_{l-1}^4+C\vv^2 h_{l-1}^2+C\vv^4h_{l-1}.
\end{split}
\end{equation*}
\end{theorem}
{\bf Proof.} By the Taylor expansion, we see
\begin{equation*}
\begin{split}
&\Psi(X_{h_l}^\vv(t_n))-\Psi(X_{h_{l-1}}^\vv(t_n))\\
=&\int_0^1[\nabla\Psi(X_{h_{l-1}}^\vv(t_n)+s(X_{h_l}^\vv(t_n)-X_{h_{l-1}}^\vv(t_n)))]{\mbox d}s\cdot(X_{h_l}^\vv(t_n)-X_{h_{l-1}}^\vv(t_n)).
\end{split}
\end{equation*}
Moreover, we have
\begin{equation}\label{var1}
\begin{split}
&{\rm Var}\left(\int_0^1[\nabla\Psi(X_{h_{l-1}}^\vv(t_n)+s(X_{h_l}^\vv(t_n)-X_{h_{l-1}}^\vv(t_n)))]{\mbox d}s\cdot(X_{h_l}^\vv(t_n)-X_{h_{l-1}}^\vv(t_n))\right)\\
\le&a\sum\limits_{i=0}^a{\rm Var}\left(\int_0^1[\nabla_i\Psi(X_{h_{l-1}}^\vv(t_n)+s(X_{h_l}^\vv(t_n)-X_{h_{l-1}}^\vv(t_n)))]{\mbox d}s\cdot[X_{h_l}^\vv(t_n)-X_{h_{l-1}}^\vv(t_n)]_i\right),
\end{split}
\end{equation}
where $\nabla_i$ is the $i$-th component of first derivatives vector and $[X_{h_l}^\vv(t_n)-X_{h_{l-1}}^\vv(t_n)]_i$ is the $i$-th component of $X_{h_l}^\vv(t_n)-X_{h_{l-1}}^\vv(t_n)$. By Lemma \ref{tech}, it is obvious to get
\begin{equation*}
\begin{split}
\E \left[\sup_{0\le t\le T}|X_{h_{l}}^\vv(t)-Z_{h_l}(t)|^2\right] \le C \vv^2,
\end{split}
\end{equation*}
and
\begin{equation*}
\begin{split}
\E \left[\sup_{0\le t\le T}|X_{h_{l-1}}^\vv(t)-Z_{h_{l-1}}(t)|^2\right] \le C \vv^2.
\end{split}
\end{equation*}
Thus, application of Lemma \ref{tech1} leads to
\begin{equation*}
\begin{split}
{\rm Var}\left(\int_0^1[\nabla_i\Psi(X_{h_{l-1}}^\vv(t_n)+s(X_{h_l}^\vv(t_n)-X_{h_{l-1}}^\vv(t_n)))]{\mbox d}s\right)\le C\vv^2.
\end{split}
\end{equation*}
Then by Theorem \ref{election} and Lemma \ref{tech2}, we get
\begin{equation}\label{var2}
\begin{split}
&{\rm Var}\left(\int_0^1[\nabla_i\Psi(X_{h_{l-1}}^\vv(t_n)+s(X_{h_l}^\vv(t_n)-X_{h_{l-1}}^\vv(t_n)))]{\mbox d}s\cdot[X_{h_l}^\vv(t_n)-X_{h_{l-1}}^\vv(t_n)]_i\right)\\
\le&C\vv^2(M^2h_l^2+\vv^4Mh_l)+C{\rm Var}\left([X_{h_l}^\vv(t_n)-X_{h_{l-1}}^\vv(t_n)]_i\right).
\end{split}
\end{equation}
Now we concentrate on ${\rm Var}\left([X_{h_l}^\vv(t_n)-X_{h_{l-1}}^\vv(t_n)]_i\right)$. For $n\le M^{l-1}-1$, $i=1,2,\cdots, a$
\begin{equation*}
\begin{split}
&[Y^\vv_{h_{l}}(t_{n+1})-Y_{h_{l-1}}^\vv(t_{n+1})]_i=[Y^\vv_{h_{l}}(t_{n})-Y_{h_{l-1}}^\vv(t_{n})]_i\\
&+h_l\sum_{k=0}^{M-1}[f_i(X^{\vv}_{h_l}(t_n^k), X^{\vv}_{h_l}(t_{n}^{k}-m_lh_l))-f_i(X^{\vv}_{h_{l}}(t_n), X^{\vv}_{h_{l}}(t_{n}-m_lh_l))]\\
&+h_l\sum_{k=0}^{M-1}[f_i(X^{\vv}_{h_l}(t_n), X^{\vv}_{h_{l}}(t_{n}-m_lh_l))-f_i(X^{\vv}_{h_{l-1}}(t_n), X^{\vv}_{h_{l-1}}(t_{n}-m_lh_l))]\\
&+\vv \sqrt{h_l}\sum_{k=0}^{M-1}[g_i(X^{\vv}_{h_l}(t_n^k), X^{\vv}_{h_l}(t_{n}^{k}-m_lh_l))-g_i(X^{\vv}_{h_{l}}(t_n), X^{\vv}_{h_{l}}(t_{n}-m_lh_l))]\xi_n^k\\
&+\vv \sqrt{h_l}\sum_{k=0}^{M-1}[g_i(X^{\vv}_{h_l}(t_n), X^{\vv}_{h_l}(t_{n}-m_lh_l))-g_i(X^{\vv}_{h_{l-1}}(t_n), X^{\vv}_{h_{l-1}}(t_{n}-m_lh_l))]\xi_n^k.
\end{split}
\end{equation*}
where $f_i$ is the $i$-th component of $f$ and $g_i$ is the $i$-th row of $g$. By computation,
\begin{equation*}
\begin{split}
&{\rm Var}[Y^\vv_{h_{l}}(t_{n+1})-Y_{h_{l-1}}^\vv(t_{n+1})]_i\le {\rm Var}[Y^\vv_{h_{l}}(t_{n})-Y_{h_{l-1}}^\vv(t_{n})]_i\\
&+4Mh_l^2\sum_{k=0}^{M-1}{\rm Var}[f_i(X^{\vv}_{h_l}(t_n^k), X^{\vv}_{h_l}(t_{n}^{k}-m_lh_l))-f_i(X^{\vv}_{h_{l}}(t_n), X^{\vv}_{h_{l}}(t_{n}-m_lh_l))]\\
&+4M^2h_l^2{\rm Var}[f_i(X^{\vv}_{h_l}(t_n), X^{\vv}_{h_{l}}(t_{n}-m_lh_l))-f_i(X^{\vv}_{h_{l-1}}(t_n), X^{\vv}_{h_{l-1}}(t_{n}-m_lh_l))]\\
&+4\vv^2h_l{\rm Var}\left\{\sum_{k=0}^{M-1}[g_i(X^{\vv}_{h_l}(t_n^k), X^{\vv}_{h_l}(t_{n}^{k}-m_lh_l))-g_i(X^{\vv}_{h_{l}}(t_n), X^{\vv}_{h_{l}}(t_{n}-m_lh_l))]\xi_n^k\right\}\\
&+4\vv^2h_l{\rm Var}\left\{\sum_{k=0}^{M-1}[g_i(X^{\vv}_{h_l}(t_n), X^{\vv}_{h_l}(t_{n}-m_lh_l))-g_i(X^{\vv}_{h_{l-1}}(t_n), X^{\vv}_{h_{l-1}}(t_{n}-m_lh_l))]\xi_n^k\right\}\\
&+2{\rm Cov}\left([Y^\vv_{h_{l}}(t_{n})-Y_{h_{l-1}}^\vv(t_{n})]_i, h_l\sum_{k=0}^{M-1}[f_i(X^{\vv}_{h_l}(t_n^k), X^{\vv}_{h_l}(t_{n}^{k}-m_lh_l))-f_i(X^{\vv}_{h_{l}}(t_n), X^{\vv}_{h_{l}}(t_{n}-m_lh_l))]\right)\\
&+2{\rm Cov}\bigg([Y^\vv_{h_{l}}(t_{n})-Y_{h_{l-1}}^\vv(t_{n})]_i, h_l\sum_{k=0}^{M-1}[f_i(X^{\vv}_{h_l}(t_n), X^{\vv}_{h_{l}}(t_{n}-m_lh_l))\\
&~~~~~~~~~~~~~~~-f_i(X^{\vv}_{h_{l-1}}(t_n), X^{\vv}_{h_{l-1}}(t_{n}-m_lh_l))]\bigg).
\end{split}
\end{equation*}
Summing both sides, for $0\le\Lambda\le M^{l-1}-1$,
\begin{equation*}
\begin{split}
&\sup\limits_{0\le n\le \Lambda+1}{\rm Var}[Y^\vv_{h_{l}}(t_{n})-Y_{h_{l-1}}^\vv(t_{n})]_i\\
\le &4Mh_l^2\sum\limits_{j=0}^{\Lambda}\sum_{k=0}^{M-1}{\rm Var}[f_i(X^{\vv}_{h_l}(t_j^k), X^{\vv}_{h_l}(t_{j}^{k}-m_lh_l))-f_i(X^{\vv}_{h_{l}}(t_j), X^{\vv}_{h_{l}}(t_{j}-m_lh_l))]\\
&+4M^2h_l^2\sum\limits_{j=0}^{\Lambda}{\rm Var}[f_i(X^{\vv}_{h_l}(t_j), X^{\vv}_{h_{l}}(t_{j}-m_lh_l))-f_i(X^{\vv}_{h_{l-1}}(t_j), X^{\vv}_{h_{l-1}}(t_{j}-m_lh_l))]\\
&+4\vv^2h_l\sum\limits_{j=0}^{\Lambda}{\rm Var}\left\{\sum_{k=0}^{M-1}[g_i(X^{\vv}_{h_l}(t_j^k), X^{\vv}_{h_l}(t_{j}^{k}-m_lh_l))-g_i(X^{\vv}_{h_{l}}(t_j), X^{\vv}_{h_{l}}(t_{j}-m_lh_l))]\xi_j^k\right\}\\
&+4\vv^2h_l\sum\limits_{j=0}^{\Lambda}{\rm Var}\left\{\sum_{k=0}^{M-1}[g_i(X^{\vv}_{h_l}(t_j), X^{\vv}_{h_l}(t_{j}-m_lh_l))-g_i(X^{\vv}_{h_{l-1}}(t_j), X^{\vv}_{h_{l-1}}(t_{j}-m_lh_l))]\xi_j^k\right\}\\
&+2\sum\limits_{j=0}^{\Lambda}{\rm Cov}\bigg([Y^\vv_{h_{l}}(t_{j})-Y_{h_{l-1}}^\vv(t_{j})]_i, h_l\sum_{k=0}^{M-1}[f_i(X^{\vv}_{h_l}(t_j^k), X^{\vv}_{h_l}(t_{j}^{k}-m_lh_l))\\
&~~~~~~~~~~~~~~~~~~-f_i(X^{\vv}_{h_{l}}(t_j), X^{\vv}_{h_{l}}(t_{j}-m_lh_l))]\bigg)\\
&+2\sum\limits_{j=0}^{\Lambda}{\rm Cov}\bigg([Y^\vv_{h_{l}}(t_{j})-Y_{h_{l-1}}^\vv(t_{j})]_i, h_l\sum_{k=0}^{M-1}[f_i(X^{\vv}_{h_l}(t_j), X^{\vv}_{h_{l}}(t_{j}-m_lh_l))\\
&~~~~~~~~~~~~~~~~~~-f_i(X^{\vv}_{h_{l-1}}(t_j), X^{\vv}_{h_{l-1}}(t_{j}-m_lh_l))]\bigg)\\
:=&I_1+I_2+I_3+I_4+I_5+I_6.
\end{split}
\end{equation*}

\begin{lemma}\label{i1}
There exists a positive constant $C$ such that
\begin{equation*}
\begin{split}
I_1\le C\vv^2M^2h_l^2+CM^5h_l^5.
\end{split}
\end{equation*}
\end{lemma}
{\bf Proof.} By the Taylor expansion,
\begin{equation*}
\begin{split}
&f_i(X^{\vv}_{h_l}(t_j^k), X^{\vv}_{h_l}(t_{j}^{k}-m_lh_l))-f_i(X^{\vv}_{h_{l}}(t_j), X^{\vv}_{h_{l}}(t_{j}-m_lh_l))\\
=& \int_0^1\big\{\nabla f_i((X^\vv_{h_l}(t_j), X^\vv_{h_l}(t_{j}-m_lh_l)) +s[(X^\vv_{h_l}(t_j^k), X^\vv_{h_l}(t_{j}^k-m_lh_l))\\
&-(X^\vv_{h_l}(t_j), X^\vv_{h_l}(t_{j}-m_lh_l))])\big\}ds\times \begin{pmatrix}
X^\vv_{h_l}(t_j^k)-X^\vv_{h_l}(t_j)\\
X^\vv_{h_l}(t_{j}^k-m_lh_l)-X^\vv_{h_l}(t_{j}-m_lh_l)
\end{pmatrix}.
\end{split}
\end{equation*}
By virtue of properties of expectation, for $r=1,2,\cdots, 2a$,
\begin{equation}\label{yu}
\begin{split}
&{\rm Var}\left[f_i(X^{\vv}_{h_l}(t_j^k), X^{\vv}_{h_l}(t_{j}^{k}-m_lh_l))-f_i(X^{\vv}_{h_{l}}(t_j), X^{\vv}_{h_{l}}(t_{j}-m_lh_l))\right]\\
\le& 2a\sum\limits_{r=1}^{2a}{\rm Var}\bigg[\int_0^1\big\{\nabla_r f_i((X^\vv_{h_l}(t_j), X^\vv_{h_l}(t_{j}-m_lh_l)) +s[(X^\vv_{h_l}(t_j^k), X^\vv_{h_l}(t_{j}^k-m_lh_l))\\
&-(X^\vv_{h_l}(t_j), X^\vv_{h_l}(t_{j}-m_lh_l))])\big\}ds\times \begin{pmatrix}
X^\vv_{h_l}(t_j^k)-X^\vv_{h_l}(t_j)\\
X^\vv_{h_l}(t_{j}^k-m_lh_l)-X^\vv_{h_l}(t_{j}-m_lh_l)
\end{pmatrix}_r\bigg],
\end{split}
\end{equation}
where $\nabla_r f_i$ is the $r$-th component of first derivatives to $f_i$, and $(\cdot)_r$ is the $r$-th component of a vector. We apply Lemma \ref{tech} and Lemma \ref{tech1} to get
\begin{equation*}
\begin{split}
&{\rm Var}\bigg[\int_0^1\big\{\nabla_r f_i((X^\vv_{h_l}(t_j), X^\vv_{h_l}(t_{j}-m_lh_l))+s[(X^\vv_{h_l}(t_j^k), X^\vv_{h_l}(t_{j}^k-m_lh_l))\\
&-(X^\vv_{h_l}(t_j), X^\vv_{h_l}(t_{j}-m_lh_l))])\big\}ds\bigg]\le C\vv^2.
\end{split}
\end{equation*}
Moreover, for $r=1,\cdots,a$, taking advantage of assumption (H) and Lemmas \ref{pmoment}-\ref{tech},
\begin{equation*}
\begin{split}
&{\rm Var}\bigg[\begin{pmatrix}
X^\vv_{h_l}(t_j^k)-X^\vv_{h_l}(t_j)\\
X^\vv_{h_l}(t_{j}^k-m_lh_l)-X^\vv_{h_l}(t_{j}-m_lh_l)
\end{pmatrix}_r\bigg]\\
\le& 3{\rm Var}\left(\sum_{q=0}^{k-1}f_r(X^{\vv}_{h_l}(t_j^q), X^{\vv}_{h_l}(t_{j}^q-m_lh_l))h_l\right)+3{\rm Var}\left(\vv \sqrt{h_l}\sum_{q=0}^{k-1}g _r(X^{\vv}_{h_l}(t_j^q), X^{\vv}_{h_l}(t_{j}^q-m_lh_l)) \xi_j^q \right)\\
&+3{\rm Var}\left[\theta f_r(X^{\vv}_{h_{l}}(t_{j}^k),X^{\vv}_{h_{l}}(t_{j}^k-m_lh_l))h_l-\theta f_r(X^{\vv}_{h_{l}}(t_{j}),X^{\vv}_{h_{l}}(t_{j}-m_lh_l))h_l\right]\\
\le& 3h_l^2{\rm Var}\left(\sum_{q=0}^{k-1}[f_r(X^{\vv}_{h_l}(t_j^q), X^{\vv}_{h_l}(t_{j}^q-m_lh_l))-f_r(Z_{h_l}(t_j^q), Z_{h_l}(t_{j}^q-m_lh_l))]\right)\\
&+3\vv^2h_l\mathbb{E}\left(\sum_{q=0}^{k-1}g_r(X^{\vv}_{h_l}(t_j^q), X^{\vv}_{h_l}(t_{j}^q-m_lh_l)) \xi_j^q \right)^2\\
&+9{\rm Var}\left[\theta f_r(X^{\vv}_{h_{l}}(t_{j}^k),X^{\vv}_{h_{l}}(t_{j}^k-m_lh_l))h_l-\theta f_r(Z_{h_{l}}(t_{j}^k),Z_{h_{l}}(t_{j}^k-m_lh_l))h_l\right]\\
&+9{\rm Var}\left[\theta f_r(Z_{h_{l}}(t_{j}^k),Z_{h_{l}}(t_{j}^k-m_lh_l))h_l-\theta f_r(Z_{h_{l}}(t_{j}),Z_{h_{l}}(t_{j}-m_lh_l))h_l\right]\\
&+9{\rm Var}\left[\theta f_r(Z_{h_{l}}(t_{j}),Z_{h_{l}}(t_{j}-m_lh_l))h_l-\theta f_r(X^{\vv}_{h_{l}}(t_{j}),X^{\vv}_{h_{l}}(t_{j}-m_lh_l))h_l\right]\\
\le&C\vv^2M^2h_l^2+C\vv^2Mh_l+C\vv^2h_l^2+Ch_l^4.
\end{split}
\end{equation*}
Similarly, for $r=a+1,\cdots, 2a$,
\begin{equation*}
\begin{split}
&{\rm Var}\bigg[\begin{pmatrix}
X^\vv_{h_l}(t_j^k)-X^\vv_{h_l}(t_j)\\
X^\vv_{h_l}(t_{j}^k-m_lh_l)-X^\vv_{h_l}(t_{j}-m_lh_l)
\end{pmatrix}_r\bigg]\\
\le&3{\rm Var}\left(
\sum_{q=0}^{k-1}f_{r-a}(X^{\vv}_{h_l}(t_{j}^q-m_lh_l), X^{\vv}_{h_l}(t_{j}^q-2m_lh_l))h_l\right)\\
&+3{\rm Var}\left(\vv \sqrt{h_l}\sum_{q=0}^{k-1}g_{r-a} (X^{\vv}_{h_l}(t_{j}^q-m_lh_l), X^{\vv}_{h_l}(t_{j}^q-2m_lh_l)) \xi_j^q\right)\\
&+3{\rm Var}\left[\theta f_{r-a}(X^{\vv}_{h_{l}}(t_{j}^k-m_lh_l),X^{\vv}_{h_{l}}(t_{j}^k-2m_lh_l))h_l-\theta f_{r-a}(X^{\vv}_{h_{l}}(t_{j}-m_lh_l),X^{\vv}_{h_{l}}(t_{j}-2m_lh_l))h_l\right]\\
\le&C\vv^2M^2h_l^2+C\vv^2Mh_l+Ch_l^4.
\end{split}
\end{equation*}
Thus, combining \eqref{yu} and Lemmas \ref{miss}, \ref{tech2}, we see
\begin{equation*}
\begin{split}
&{\rm Var}\left[f_i(X^{\vv}_{h_l}(t_j^k), X^{\vv}_{h_l}(t_{j}^{k}-m_lh_l))-f_i(X^{\vv}_{h_{l}}(t_j), X^{\vv}_{h_{l}}(t_{j}-m_lh_l))\right]\\
\le& 2a\sum\limits_{r=1}^{2a}\bigg\{C\vv^2(M^2h_l^2+\vv^2Mh_l)+C{\rm Var}\bigg[\begin{pmatrix}
X^\vv_{h_l}(t_j^k)-X^\vv_{h_l}(t_j)\\
X^\vv_{h_l}(t_{j}^k-m_lh_l)-X^\vv_{h_l}(t_{j}-m_lh_l)
\end{pmatrix}_r\bigg]\bigg\}\\
\le&C\vv^2Mh_l+Ch_l^4,
\end{split}
\end{equation*}
which leads to
\begin{equation*}
\begin{split}
I_1\le C\vv^2M^2h_l^2+CMh_l^5.
\end{split}
\end{equation*}\hfill $\Box$

\begin{lemma}\label{i2}
There exists a positive constant $C$ such that
\begin{equation*}
\begin{split}
I_2\le& C\vv^2M^3h_l^3+C\vv^6M^2h_l^2+CM^2h_l^2\sum\limits_{j=0}^{\Lambda}\sum\limits_{i=0}^{a}\sup\limits_{0\le n\le j}{\rm Var}\left([X_{h_l}^\vv(t_n)-X_{h_{l-1}}^\vv(t_n)]_i\right).
\end{split}
\end{equation*}
\end{lemma}
{\bf Proof.} Application of the Taylor expansion gives that
\begin{equation*}
\begin{split}
&f_i(X^{\vv}_{h_l}(t_j), X^{\vv}_{h_{l}}(t_{j}-m_lh_l))-f_i(X^{\vv}_{h_{l-1}}(t_j), X^{\vv}_{h_{l-1}}(t_{j}-m_lh_l))\\
=& \int_0^1\big\{\nabla f_i((X^\vv_{h_{l-1}}(t_j), X^\vv_{h_{l-1}}(t_{j}-m_lh_l)) +s[(X^\vv_{h_l}(t_j), X^\vv_{h_l}(t_{j}-m_lh_l))\\
&-(X^\vv_{h_{l-1}}(t_j), X^\vv_{h_{l-1}}(t_{j}-m_lh_l))])\big\}ds\times \begin{pmatrix}
X^\vv_{h_l}(t_j)-X^\vv_{h_{l-1}}(t_j)\\
X^\vv_{h_l}(t_{j}-m_lh_l)-X^\vv_{h_{l-1}}(t_{j}-m_lh_l)
\end{pmatrix}.
\end{split}
\end{equation*}
Taking similar steps as in Lemma \ref{i1}, Lemma \ref{tech} together with Lemma \ref{tech1} yield
\begin{equation*}
\begin{split}
&{\rm Var}\bigg[\int_0^1\big\{\nabla_r f_i((X^\vv_{h_{l-1}}(t_j), X^\vv_{h_{l-1}}(t_{j}-m_lh_l)) +s[(X^\vv_{h_l}(t_j), X^\vv_{h_l}(t_{j}-m_lh_l))\\
&-(X^\vv_{h_{l-1}}(t_j), X^\vv_{h_{l-1}}(t_{j}-m_lh_l))])\big\}ds\bigg]\le C\vv^2
\end{split}
\end{equation*}
for $r=1,2,\cdots, 2a$. Further, Theorem \ref{election} and Lemma \ref{tech2} yield
\begin{equation*}
\begin{split}
I_2\le& CMh_l\vv^2(M^2h_l^2+\vv^4Mh_l)+CM^2h_l^2\sum\limits_{j=0}^{\Lambda}\sum\limits_{r=1}^{2a}{\rm Var}\bigg[\begin{pmatrix}
X^\vv_{h_l}(t_j)-X^\vv_{h_{l-1}}(t_j)\\
X^\vv_{h_l}(t_{j}-m_lh_l)-X^\vv_{h_{l-1}}(t_{j}-m_lh_l)
\end{pmatrix}_r\bigg]\\
\le& C\vv^2 M^3h_l^3+C\vv^6M^2h_l^2+CM^2h_l^2\sum\limits_{j=0}^{\Lambda}\sum\limits_{i=0}^{a}\sup\limits_{0\le n\le j}{\rm Var}\left([X_{h_l}^\vv(t_n)-X_{h_{l-1}}^\vv(t_n)]_i\right).
\end{split}
\end{equation*}\hfill $\Box$

\begin{lemma}\label{i3}
There exists a positive constant $C$ such that
\begin{equation*}
\begin{split}
I_3\le C\vv^2M^2h_l^2+C\vv^4Mh_l.
\end{split}
\end{equation*}
\end{lemma}
{\bf Proof.} By assumption (H) and Lemma \ref{miss},
\begin{equation*}
\begin{split}
I_3\le& 4\vv^2h_l\sum\limits_{j=0}^{\Lambda}\mathbb{E}\left|\sum_{k=0}^{M-1}[g_i(X^{\vv}_{h_l}(t_j^k), X^{\vv}_{h_l}(t_{j}^{k}-m_lh_l))-g_i(X^{\vv}_{h_{l}}(t_j), X^{\vv}_{h_{l}}(t_{j}-m_lh_l))]\xi_j^k\right|^2\\
\le&C\vv^2(M^2h_l^2+\vv^2Mh_l).
\end{split}
\end{equation*}
This completes the proof.\hfill $\Box$

\begin{lemma}\label{i4}
There exists a positive constant $C$ such that
\begin{equation*}
\begin{split}
I_4\le C\vv^2M^2h_l^2+C\vv^6Mh_l.
\end{split}
\end{equation*}
\end{lemma}
{\bf Proof.} By assumption (H) and Theorem \ref{election},
\begin{equation*}
\begin{split}
I_4\le& 4\vv^2h_l\sum\limits_{j=0}^{\Lambda}\mathbb{E}\left|\sum_{k=0}^{M-1}[g_i(X^{\vv}_{h_l}(t_j), X^{\vv}_{h_l}(t_{j}-m_lh_l))-g_i(X^{\vv}_{h_{l-1}}(t_j), X^{\vv}_{h_{l-1}}(t_{j}-m_lh_l))]\xi_j^k\right|^2\\
\le&C\vv^2(M^2h_l^2+\vv^4Mh_l).
\end{split}
\end{equation*}\hfill $\Box$

\begin{lemma}\label{i5}
There exists a positive constant $C$ such that
\begin{equation*}
\begin{split}
I_5\le C\vv^2M^2h_l^2+CM^4h_l^4+\frac{1}{4}\sup_{0\le n\le \Lambda+1}{\rm Var}\left([Y_{h_l}^\vv(t_n)-Y_{h_{l-1}}^\vv(t_n)]_i\right).
\end{split}
\end{equation*}
\end{lemma}
{\bf Proof.} Recall from Lemma \ref{taylor} that
\begin{equation}\label{lyon}
\begin{split}
I_5=&2\sum\limits_{j=0}^{\Lambda}{\rm Cov}\left([Y^\vv_{h_{l}}(t_{j})-Y_{h_{l-1}}^\vv(t_{j})]_i, h_l\sum_{k=0}^{M-1}(A_j^{ik}+B_j^{ik}+C_j^{ik})\right)\\
=&2\sum\limits_{j=0}^{\Lambda}{\rm Cov}\left([Y^\vv_{h_{l}}(t_{j})-Y_{h_{l-1}}^\vv(t_{j})]_i, h_l\sum_{k=0}^{M-1}(A_j^{ik}+C_j^{ik})\right)\\
\le&\frac{1}{4}\sup_{0\le n\le \Lambda+1}{\rm Var}\left([Y_{h_l}^\vv(t_n)-Y_{h_{l-1}}^\vv(t_n)]_i\right)+Ch_l\sum\limits_{j=0}^{\Lambda}\sum_{k=0}^{M-1}{\rm Var}(A_j^{ik})+Ch_l\sum\limits_{j=0}^{\Lambda}\sum_{k=0}^{M-1}{\rm Var}(C_j^{ik}).
\end{split}
\end{equation}
We are now going to estimate ${\rm Var}(A_j^{ik})$ and ${\rm Var}(C_j^{ik})$. By the Taylor expansion,
\begin{equation*}
\begin{split}
{\rm Var}(A_j^{ik})=&{\rm Var}\bigg[\int_0^1\big\{\nabla f_i((X^\vv_{h_l}(t_j), X^\vv_{h_l}(t_{j}-m_lh_l)) +s[(X^\vv_{h_l}(t_j^k), X^\vv_{h_l}(t_{j}^k-m_lh_l))\\
&-(X^\vv_{h_l}(t_j), X^\vv_{h_l}(t_{j}-m_lh_l))])\big\}ds\times\begin{pmatrix}\sigma_j^{11}+\sigma_j^{13}\\
\sigma_j^{21}+\sigma_j^{23}
\end{pmatrix}\bigg]\\
\le&2a\sum\limits_{r=1}^{2a}{\rm Var}\bigg[\int_0^1\big\{\nabla_r f_i((X^\vv_{h_l}(t_j), X^\vv_{h_l}(t_{j}-m_lh_l)) +s[(X^\vv_{h_l}(t_j^k), X^\vv_{h_l}(t_{j}^k-m_lh_l))\\
&-(X^\vv_{h_l}(t_j), X^\vv_{h_l}(t_{j}-m_lh_l))])\big\}ds\times\begin{pmatrix}\sigma_j^{11}+\sigma_j^{13}\\
\sigma_j^{21}+\sigma_j^{23}
\end{pmatrix}_r\bigg]\\
\le&C\vv^2(M^2h_l^2+\vv^2M^3h_l^3)+C\sum\limits_{r=1}^{2a}{\rm Var}\bigg[\begin{pmatrix}\sigma_j^{11}+\sigma_j^{13}\\
\sigma_j^{21}+\sigma_j^{23}
\end{pmatrix}_r\bigg]\\
\le&C\vv^2M^2h_l^2+C\vv^4M^3h_l^3+Ch_l^4.
\end{split}
\end{equation*}
Since for $r=1,\cdots,a$, similar to the procedure of Lemma \ref{i1}, by assumption (H) and Lemmas \ref{pmoment}-\ref{tech}, we get
\begin{equation*}
\begin{split}
&{\rm Var}\bigg[\begin{pmatrix}
\sigma_j^{11}+\sigma_j^{13}\\
\sigma_j^{21}+\sigma_j^{23}
\end{pmatrix}_r\bigg]\le 2{\rm Var}\left(\sum_{q=0}^{k-1}f_r(X^{\vv}_{h_l}(t_j^q), X^{\vv}_{h_l}(t_{j}^q-m_lh_l))h_l\right)\\
&+2{\rm Var}\left[\theta f_r(X^{\vv}_{h_{l}}(t_{j}^k),X^{\vv}_{h_{l}}(t_{j}^k-m_lh_l))h_l-\theta f_r(X^{\vv}_{h_{l}}(t_{j}),X^{\vv}_{h_{l}}(t_{j}-m_lh_l))h_l\right]\\
\le&C\vv^2M^2h_l^2+C\vv^2h_l^2+Ch_l^4.
\end{split}
\end{equation*}
Similarly, for $r=a+1,\cdots, 2a$,
\begin{equation*}
\begin{split}
&{\rm Var}\bigg[\begin{pmatrix}
\sigma_j^{11}+\sigma_j^{13}\\
\sigma_j^{21}+\sigma_j^{23}
\end{pmatrix}_r\bigg]\le C\vv^2M^2h_l^2+C\vv^2h_l^2+Ch_l^4.
\end{split}
\end{equation*}
Similar to the estimation of $\mathbb{E}|C_j^{k}|^2$ in Theorem \ref{election}, we easily get
\begin{equation*}
\begin{split}
{\rm Var}(C_j^{ik})\le\mathbb{E}|C_j^{ik}|^2\le C\vv^2M^3h_l^3+C\vv^4M^2h_l^2.
\end{split}
\end{equation*}
Then, we derive from \eqref{lyon} that
\begin{equation*}
\begin{split}
I_5\le&\frac{1}{4}\sup_{0\le n\le \Lambda+1}{\rm Var}\left([Y_{h_l}^\vv(t_n)-Y_{h_{l-1}}^\vv(t_n)]_i\right)\\
&+C(\vv^2M^2h_l^2+\vv^4M^3h_l^3+h_l^4)+C(\vv^2M^3h_l^3+\vv^4M^2h_l^2).
\end{split}
\end{equation*}\hfill $\Box$

\begin{lemma}\label{i6}
There exists a positive constant $C$ such that
\begin{equation*}
\begin{split}
I_6\le& C\vv^2 M^2h_l^2+C\vv^6Mh_l+\frac{1}{4}\sup_{0\le n\le \Lambda+1}{\rm Var}\left([Y_{h_l}^\vv(t_n)-Y_{h_{l-1}}^\vv(t_n)]_i\right)\\
&+CMh_l\sum\limits_{j=0}^\Lambda\sum\limits_{i=1}^{a}\sup\limits_{0\le n\le j}{\rm Var}\left([X_{h_l}^\vv(t_n)-X_{h_{l-1}}^\vv(t_n)]_i\right).
\end{split}
\end{equation*}
\end{lemma}
{\bf Proof.} Obviously, by the result of Lemma \ref{i2}
\begin{equation*}
\begin{split}
I_6\le& \frac{1}{4}\sup_{0\le n\le \Lambda+1}{\rm Var}\left([Y_{h_l}^\vv(t_n)-Y_{h_{l-1}}^\vv(t_n)]_i\right)\\
&+CMh_l\sum\limits_{j=0}^\Lambda{\rm Var}\left(f_i(X^{\vv}_{h_l}(t_j), X^{\vv}_{h_{l}}(t_{j}-m_lh_l))-f_i(X^{\vv}_{h_{l-1}}(t_j), X^{\vv}_{h_{l-1}}(t_{j}-m_lh_l))\right)\\
\le&\frac{1}{4}\sup_{0\le n\le \Lambda+1}{\rm Var}\left([Y_{h_l}^\vv(t_n)-Y_{h_{l-1}}^\vv(t_n)]_i\right)+C(\vv^2 M^2h_l^2+C\vv^6Mh_l)\\
&+CMh_l\sum\limits_{j=0}^\Lambda\sum\limits_{i=1}^{a}\sup\limits_{0\le n\le j}{\rm Var}\left([X_{h_l}^\vv(t_n)-X_{h_{l-1}}^\vv(t_n)]_i\right).
\end{split}
\end{equation*}\hfill $\Box$\\
{\bf Continue of Theorem \ref{th1}.} By Lemmas \ref{i1}-\ref{i6}, we see
\begin{equation}\label{wind}
\begin{split}
&\sup\limits_{0\le n\le \Lambda+1}{\rm Var}[Y^\vv_{h_{l}}(t_{n})-Y_{h_{l-1}}^\vv(t_{n})]_i\le C\vv^2M^2h_l^2+C\vv^4Mh_l+CM^4h_l^4\\
&+CMh_l\sum\limits_{j=0}^\Lambda\sum\limits_{i=1}^{a}\sup\limits_{0\le n\le j}{\rm Var}\left([X_{h_l}^\vv(t_n)-X_{h_{l-1}}^\vv(t_n)]_i\right).
\end{split}
\end{equation}
Since by Lemma \ref{tech}, we have
\begin{equation}\label{home}
\begin{split}
&{\rm Var}[X^\vv_{h_{l}}(t_{n})-X_{h_{l-1}}^\vv(t_{n})]_i\\
\le&3{\rm Var}[Y^\vv_{h_{l}}(t_{n})-Y_{h_{l-1}}^\vv(t_{n})]_i+3\theta^2h_l^2{\rm Var}[f_i(X_{h_l}^\vv(t_{n}),X_{h_l}^\vv(t_{n}-m_lh_l))]\\
&+3\theta^2h_{l-1}^2{\rm Var}[f_i(X_{h_{l-1}}^\vv(t_{n}),X_{h_{l-1}}^\vv(t_{n}-m_lh_l))]\\
\le&3{\rm Var}[Y^\vv_{h_{l}}(t_{n})-Y_{h_{l-1}}^\vv(t_{n})]_i\\
&+3\theta^2h_l^2{\rm Var}[f_i(X_{h_l}^\vv(t_{n}),X_{h_l}^\vv(t_{n}-m_lh_l))-f_i(Z_{h_l}(t_{n}),Z_{h_l}(t_{n}-m_lh_l))]\\
&+3\theta^2h_{l-1}^2{\rm Var}[f_i(X_{h_{l-1}}^\vv(t_{n}),X_{h_{l-1}}^\vv(t_{n}-m_lh_l))-f_i(Z_{h_{l-1}}(t_{n}),Z_{h_{l-1}}(t_{n}-m_lh_l))]\\
\le&3{\rm Var}[Y^\vv_{h_{l}}(t_{n})-Y_{h_{l-1}}^\vv(t_{n})]_i+C\vv^2M^2h_l^2.
\end{split}
\end{equation}
Then, by \eqref{wind} and \eqref{home}, for $\Lambda\le M^{l-1}-1$
\begin{equation*}
\begin{split}
&\sup\limits_{0\le n\le \Lambda+1}{\rm Var}[X^\vv_{h_{l}}(t_{n})-X_{h_{l-1}}^\vv(t_{n})]_i\\
\le& C\vv^2M^2h_l^2+C\vv^4Mh_l+CM^4h_l^4+CMh_l\sum\limits_{j=0}^\Lambda\sum\limits_{i=1}^{a}\sup\limits_{0\le n\le j}{\rm Var}\left([X_{h_l}^\vv(t_n)-X_{h_{l-1}}^\vv(t_n)]_i\right).
\end{split}
\end{equation*}
The Gronwall inequality leads to
\begin{equation*}
\begin{split}
&\sup\limits_{0\le n<M^{l-1}}\sup\limits_{0\le i\le a}{\rm Var}[X^\vv_{h_{l}}(t_{n})-X_{h_{l-1}}^\vv(t_{n})]_i\le C\vv^2M^2h_l^2+C\vv^4Mh_l+CM^4h_l^4.
\end{split}
\end{equation*}
The desired result then follows from \eqref{var1}-\eqref{var2}. \hfill $\Box$

\begin{remark}
{\rm From Section \ref{sec2.1}, we see that the theta EM scheme has the following property
\begin{equation*}
\sup_{0\le n<M^{l-1}}{\rm Var}(\Psi(X_{h_l}^\vv(t_n))-\Psi(X_{h_{l-1}}^\vv(t_n)))\le Ch_{l-1}^2+C\vv^2h_{l-1},
\end{equation*}
while for the multilevel Monte Carlo theta EM scheme, the variance is bounded by $\mathcal{O}(h_{l-1}^4+\vv^2h_{l-1}^2+\vv^4h_{l-1})$. That is, the multilevel Monte Carlo theta EM scheme is more efficient than the theta EM scheme.
}
\end{remark}

\section{SDDEs under One-side Lipschitz Condition}
In this section, instead of the global Lipschitz condition (H), we impose weaker assumptions to \eqref{2.0}. We assume that:
 \begin{description}
\item[(H1)] There exist $\alpha_1, \alpha_2>1$ such that for some $p\ge 2$, $r\ge1$
\begin{displaymath}
2\langle x-\bar{x}, f(x,y)-f (\bar{x},\bar{y})\rangle+(p-1)\vv^2|g(x,y)-g (\bar{x},\bar{y})|^2\le \alpha_1(|x-\bar{x}|^2+|y-\bar{y}|^2)
\end{displaymath}
and
\begin{displaymath}
|f(x,y)-f (\bar{x},\bar{y})|\le \alpha_2(1+|x|^r+|\bar{x}|^r+|y|^r+|\bar{y}|^r)(|x-\bar{x}|+|y-\bar{y}|)
\end{displaymath}
for all $x, y, \bar x, \bar y \in \RR^a$.
\item[(H2)] There exists a positive constant $\alpha_3$ such that
\begin{displaymath}
|g(x,y)|^2\le \alpha_3(1+|x|^2+|y|^2)
\end{displaymath}
for all $x, y\in \RR^a$.
\end{description}

\begin{lemma}
Let assumptions (H1) and (H2) hold. Then, for any $T>0$ and $p\ge 2,$ we have
$$
\sup_{0\le t\le T}\E|X^\vv(t)|^p\le C.
$$
\end{lemma}

\begin{remark}%\label{onem}
{\rm Assumption (H1) implies that for any $x, y\in \RR^a$
\begin{equation}\label{mono}
\langle x,f(x,y)\rangle\le \bar{\alpha}_1(1+|x|^2+|y|^2)
\end{equation}
where $\bar{\alpha}_1=\alpha_1\vee\frac{1}{2}|f(0,0)|^2$. Moreover, assumption (H1) also implies that
\begin{equation}\label{dec}
\begin{split}
&\quad (p-1)\vv^2|g(x, y)-g(\bar{x}, \bar{y})|^2\\
&\le \alpha_1(|x-\bar{x}|^2+|y-\bar{y}|^2)+2|x-\bar{x}||f(x, y)-f(\bar{x}, \bar{y})|\\
&\le \tilde{\alpha}(1+|x|^r+|\bar{x}|^r+|y|^r+|\bar{y}|^r)(|x-\bar{x}|^2+|y-\bar{y}|^2)
\end{split}
\end{equation}
where $\tilde{\alpha}$ is a constant depends on $\alpha_1$ and $\alpha_2$.
}
\end{remark}

\begin{remark}%\label{twom}
{\rm Assumptions (H1)-(H2) guarantee the existence and uniqueness of the solution to \eqref{2.0}.
}
\end{remark}

In order to guarantee the finiteness of $p$-th moment of the numerical solutions to \eqref{2.0}, we make a tiny modification to the drift coefficient. Given any $T>0$, let $M\ge 2, l>1$, $h_l=T\cdot M^{-l}, h_{l-1}=T\cdot M^{-(l-1)},$ define
\begin{equation}\label{delta}
\begin{split}
&f_{h_l}(x,y):=\frac{f(x,y)}{1+{h_{l-1}}^{\delta}|f(x,y)|}
\end{split}
\end{equation}
for any $x, y\in \RR^a$ and some $\delta\in(0,\frac{1}{2}]$.

\begin{remark}
{\rm With the definition of $f_{h_{l}}$, it is easy to show that under assumption (H1) the following condition hold:
\begin{equation}\label{delbound}
\begin{split}
&|f_{h_l}(x,y)|\le \min\left(|f(x,y)|,h_{l-1}^{-\delta}\right).
\end{split}
\end{equation}
Moreover, one can verify that $f_{h_{l}}$ satisfies the following properties:
\begin{equation}\label{delones}
\left\langle x-\bar{x}, f_{h_l}(x,y)-f_{h_l}(\bar{x},\bar{y})\right\rangle\le \frac{\alpha_1}{2}(|x-\bar{x}|^2+|y-\bar{y}|^2),
\end{equation}
and
\begin{equation}\label{delmono}
\langle x,f_{h_{l}}(x,y)\rangle \le \bar{\alpha}_1(1+|x|^2+|y|^2).
\end{equation}
Furthermore,
\begin{equation}\label{del}
|f(x,y)-f_{h_l}(x,y)|^p\le \bar{\alpha}_2 h_{l-1}^{\delta p}\big[1+|x|^{2(r+1)p}+|y|^{2(r+1)p}\big]
\end{equation}
where $\bar{\alpha}_2=[\alpha_2+|f(0,0)|]^{2p}$.}
\end{remark}
Similar to the global Lipschitz case, assume there exists an $m_l$ such that $\tau=m_lh_l$ and define
\begin{equation}\label{1c1}
\begin{split}
&X^{\vv}_{h_l}(t)-\theta f_{h_l}(X^{\vv}_{h_l}(t),X^{\vv}_{h_l}(t-\tau))h_l\\
=&\xi(0)-\theta f_{h_l}(\xi(0),\xi(-\tau))h_l+\int_0^tf_{h_l}(X^{\vv}_{h_l}(\eta_{h_{l}}(s)), X^{\vv}_{h_l}(\eta_{h_{l}}(s-\tau))){\mbox d}s\\
&+\vv\int_0^tg(X^{\vv}_{h_l}(\eta_{h_{l}}(s)), X^{\vv}_{h_l}(\eta_{h_{l}}(s-\tau))){\mbox d}W(s),
\end{split}
\end{equation}
and
\begin{equation}\label{1c2}
\begin{split}
&X^{\vv}_{h_{l-1}}(t)-\theta f_{h_{l-1}}(X^{\vv}_{h_{l-1}}(t),X^{\vv}_{h_{l-1}}(t-\tau))h_{l-1}\\
=&\xi(0)-\theta f_{h_{l-1}}(\xi(0),\xi(-\tau))h_{l-1}+\int_0^tf_{h_{l-1}}(X^{\vv}_{h_{l-1}}(\eta_{h_{l-1}}(s)), X^{\vv}_{h_{l-1}}(\eta_{h_{l-1}}(s-\tau))){\mbox d}s\\
&+\vv\int_0^tg(X^{\vv}_{h_{l-1}}(\eta_{h_{l-1}}(s)), X^{\vv}_{h_{l-1}}(\eta_{h_{l-1}}(s-\tau))){\mbox d}W(s),
\end{split}
\end{equation}
where $\eta_{h_l}(s)=\left\lfloor s/h_l\right\rfloor h_l$ and $\eta_{h_{l-1}}(s)=\left\lfloor s/h_{l-1}\right\rfloor h_{l-1}$. Here $\theta\in[0,1]$ is a parameter to control the implicitness. For $n\in \{0, 1, \ldots, M^{l-1}-1\}$ and $k\in \{0, \ldots, M\}$, let
$$
t_n=nh_{l-1} \mbox{ and } t_n^k=nh_{l-1}+kh_l.
$$
This means we divide the interval $[t_n, t_{n+1}]$ into $M$ equal parts, we have $t_n^0=t_n, t_n^{M}=t_{n+1}.$ We can rewrite \eqref{1c1} and \eqref{1c2} as the following discretization schemes.
For  $n\in \{0, 1, \ldots, M^{l-1}-1\}$ and $k\in \{0, \ldots, M-1\}$, let
\begin{equation}\label{1c3}
\begin{split}
&X^{\vv}_{h_l}(t_n^{k+1})-\theta f_{h_{l}}(X^{\vv}_{h_{l}}(t_n^{k+1}),X^{\vv}_{h_{l}}(t_{n}^{k+1}-m_lh_l))h_l\\
=& X^{\vv}_{h_l}(t_n^{k})-\theta f_{h_{l}}(X^{\vv}_{h_{l}}(t_n^{k}),X^{\vv}_{h_{l}}(t_{n}^{k}-m_lh_l))h_l+f_{h_{l}}(X^{\vv}_{h_l}(t_n^k), X^{\vv}_{h_l}(t_{n}^{k}-m_lh_l))h_l\\
&+\vv \sqrt{h_l}g (X^{\vv}_{h_l}(t_n^k), X^{\vv}_{h_l}(t_{n}^{k}-m_lh_l))\xi_n^k,
\end{split}
\end{equation}
where the random vector $\xi_n^k\in\mathbb{R}^d$ has independent components, and each component is distributed as $N(0, 1).$ This implies
\begin{equation}\label{1c4}
\begin{split}
&X^{\vv}_{h_l}(t_{n+1})-\theta f_{h_{l}}(X^{\vv}_{h_{l}}(t_{n+1}),X^{\vv}_{h_{l}}(t_{n+1}-m_lh_l))h_l\\
=&X^{\vv}_{h_l}(t_n)-\theta f_{h_{l}}(X^{\vv}_{h_{l}}(t_{n}),X^{\vv}_{h_{l}}(t_{n}-m_lh_l))h_l+\sum_{k=0}^{M-1}f_{h_{l}}(X^{\vv}_{h_l}(t_n^k), X^{\vv}_{h_l}(t_{n}^{k}-m_lh_l))h_l\\
&+\vv \sqrt{h_l}\sum_{k=0}^{M-1}g (X^{\vv}_{h_l}(t_n^k), X^{\vv}_{h_l}(t_{n}^{k}-m_lh_l))\xi_n^k.
\end{split}
\end{equation}
To simulate $X_{h_{l-1}}^\vv,$ we use
\begin{equation}\label{1c5}
\begin{split}
&X^{\vv}_{h_{l-1}}(t_{n+1})-\theta f_{h_{l-1}}(X^{\vv}_{h_{l-1}}(t_{n+1}),X^{\vv}_{h_{l-1}}(t_{n+1}-m_lh_l))h_{l-1}\\
=&X^{\vv}_{h_{l-1}}(t_n)-\theta f_{h_{l-1}}(X^{\vv}_{h_{l-1}}(t_{n}),X^{\vv}_{h_{l-1}}(t_{n}-m_lh_l))h_{l-1}+f_{h_{l-1}}(X^{\vv}_{h_{l-1}}(t_n), X^{\vv}_{h_{l-1}}(t_{n}-m_lh_l))h_{l-1}\\
&+\vv \sqrt{h_{l}}g (X^{\vv}_{h_{l-1}}(t_n), X^{\vv}_{h_{l-1}}(t_{n}-m_lh_l))\sum_{k=0}^{M-1}\xi_n^k.
\end{split}
\end{equation}
For convenience, let
\begin{equation*}
Y^{\vv}_{h_l}(t):=X^{\vv}_{h_l}(t)-\theta f_{h_{l}}(X^{\vv}_{h_{l}}(t),X^{\vv}_{h_{l}}(t-\tau))h_l,
\end{equation*}
 and
\begin{equation*}
 Y^{\vv}_{h_{l-1}}(t):=X^{\vv}_{h_{l-1}}(t)-\theta f_{h_{l-1}}(X^{\vv}_{h_{l-1}}(t),X^{\vv}_{h_{l-1}}(t-\tau))h_{l-1}.
\end{equation*}
Furthermore, in order to ensure the existence and uniqueness of solutions to implicit equations \eqref{1c1} and \eqref{1c2}, we assume that $h_{l-1}\theta<\frac{2}{\alpha_1}$ according to the monotone operator \cite{zei}. Thus, in this section, we set $h^*\in\left(0,\frac{2}{\theta\alpha_1}\right)$, and let $h_{l-1}\in(0,h^*]$ for $\theta\in(0,1]$, while for $\theta=0$, let $h_{l-1}\in(0,1)$.

We now have the following estimates.

\begin{lemma}\label{1pmoment}
Let assumptions (H1) and (H2) hold. Then, for any $T>0$ and $p\ge 2,$ we have
$$
\sup_{0\le t\le T}\E|X_{h_{l}}^\vv(t)|^p\le C,
$$
and
$$
\sup_{0\le t\le T}\E |X_{h_{l-1}}^\vv(t)|^p\le C.
$$
\end{lemma}
{\bf Proof.} Here we concentrate on the first part, since the second part can be proved similarly. For $x>0$, let $\lfloor x\rfloor$ be the integer part of $x$. For any $t\in[0,T]$, applying the It\^{o} formula to $[1+|Y^{\vv}_{h_l}(t)|^{2}]^{\frac{p}{2}}$, we obtain
\begin{equation*}
\begin{split}
&\mathbb{E}[1+|Y^{\vv}_{h_l}(t)|^{2}]^{\frac{p}{2}}\le\mathbb{E}[1+|Y^{\vv}_{h_l}(0)|^{2}]^{\frac{p}{2}}\\
&+p\mathbb{E}\int_0^t[1+|Y^{\vv}_{h_l}(s)|^{2}]^{\frac{p-2}{2}}\langle Y^{\vv}_{h_l}(s), f_{h_l}(X^{\vv}_{h_l}(\eta_{h_{l}}(s)), X^{\vv}_{h_l}(\eta_{h_{l}}(s-\tau)))\rangle\mbox{d}s\\
&+\frac{1}{2}p(p-1)\mathbb{E}\int_0^t[1+|Y^{\vv}_{h_l}(s)|^{2}]^{\frac{p-2}{2}}|\vv g(X^{\vv}_{h_l}(\eta_{h_{l}}(s)), X^{\vv}_{h_l}(\eta_{h_{l}}(s-\tau)))|^2\mbox{d}s\\
\le&\mathbb{E}[1+|Y^{\vv}_{h_l}(0)|^{2}]^{\frac{p}{2}}+\frac{1}{2}p(p-1)\mathbb{E}\int_0^t[1+|Y^{\vv}_{h_l}(s)|^{2}]^{\frac{p-2}{2}}|\vv g(X^{\vv}_{h_l}(\eta_{h_{l}}(s)), X^{\vv}_{h_l}(\eta_{h_{l}}(s-\tau)))|^2\mbox{d}s\\
&+p\mathbb{E}\int_0^t[1+|Y^{\vv}_{h_l}(s)|^{2}]^{\frac{p-2}{2}}\langle X^{\vv}_{h_l}(\eta_{h_{l}}(s)), f_{h_l}(X^{\vv}_{h_l}(\eta_{h_{l}}(s)), X^{\vv}_{h_l}(\eta_{h_{l}}(s-\tau)))\rangle\mbox{d}s\\
&+p\mathbb{E}\int_0^t[1+|Y^{\vv}_{h_l}(s)|^{2}]^{\frac{p-2}{2}}\langle Y^{\vv}_{h_l}(s)-X^{\vv}_{h_l}(\eta_{h_{l}}(s)), f_{h_l}(X^{\vv}_{h_l}(\eta_{h_{l}}(s)), X^{\vv}_{h_l}(\eta_{h_{l}}(s-\tau)))\rangle\mbox{d}s\\
=: &\mathbb{E}[1+|Y^{\vv}_{h_l}(0)|^{2}]^{\frac{p}{2}}+E_1(t)+E_2(t)+E_3(t),
\end{split}
\end{equation*}
where $Y^{\vv}_{h_l}(0)=\xi(0)-\theta f_{h_l}(\xi(0),\xi(-\tau))h_l$. With (H2), \eqref{delbound}, \eqref{delmono} and the Young inequality, we have
\begin{equation*}
\begin{split}
E_1(t)+E_2(t)\le &C\mathbb{E}\int_0^t[1+|Y^{\vv}_{h_l}(s)|^{2}]^{\frac{p-2}{2}}\bigg(1+|X^{\vv}_{h_l}(\eta_{h_{l}}(s))|^2+|X^{\vv}_{h_l}(\eta_{h_{l}}(s-\tau))|^2\bigg)\mbox{d}s\\
\le &C+C \mathbb{E}\int_0^t\left([1+|Y^{\vv}_{h_l}(s)|^{2}]^{\frac{p}{2}}+|X^{\vv}_{h_l}(\eta_{h_{l}}(s))|^p+|X^{\vv}_{h_l}(\eta_{h_{l}}(s-\tau))|^p\right)\mbox{d}s\\
\le &C+C \mathbb{E}\int_0^t\bigg[1+|X^{\vv}_{h_l}(s)|^p+|\theta f_{h_l}(X^{\vv}_{h_l}(s),X^{\vv}_{h_l}(s-\tau))h_l|^p\\
&~~~~~~~~~~~~~~~~~~+|X^{\vv}_{h_l}(\eta_{h_{l}}(s))|^p+|X^{\vv}_{h_l}(\eta_{h_{l}}(s-\tau))|^p\bigg]\mbox{d}s\\
\le&C+Ch_{l-1}^{(1-\delta)p}+C\mathbb{E}\int_0^t\bigg(|X^{\vv}_{h_l}(s)|^p+|X^{\vv}_{h_l}(\eta_{h_{l}}(s))|^p+|X^{\vv}_{h_l}(\eta_{h_{l}}(s-\tau))|^p\bigg)\mbox{d}s\\
\le&C+C\int_0^t\sup\limits_{0\le u\le s}\mathbb{E}|X^{\vv}_{h_l}(u)|^p\mbox{d}s.
\end{split}
\end{equation*}
Furthermore, it is easy to observe that,
\begin{equation*}
\begin{split}
&E_3(t)\le p\mathbb{E}\int_0^t[1+|Y^{\vv}_{h_l}(s)|^{2}]^{\frac{p-2}{2}}\langle Y^{\vv}_{h_l}(s)-Y^{\vv}_{h_l}(\eta_{h_{l}}(s)), f_{h_l}(X^{\vv}_{h_l}(\eta_{h_{l}}(s)), X^{\vv}_{h_l}(\eta_{h_{l}}(s-\tau)))\rangle\mbox{d}s\\
=&p\mathbb{E}\int_0^t[1+|Y^{\vv}_{h_l}(\eta_{h_{l}}(s))|^{2}]^{\frac{p-2}{2}}\langle Y^{\vv}_{h_l}(s)-Y^{\vv}_{h_l}(\eta_{h_{l}}(s)), f_{h_l}(X^{\vv}_{h_l}(\eta_{h_{l}}(s)), X^{\vv}_{h_l}(\eta_{h_{l}}(s-\tau)))\rangle\mbox{d}s\\
&+p\mathbb{E}\int_0^t\left\{[1+|Y^{\vv}_{h_l}(s)|^{2}]^{\frac{p-2}{2}}-[1+|Y^{\vv}_{h_l}(\eta_{h_{l}}(s))|^{2}]^{\frac{p-2}{2}}\right\}\langle Y^{\vv}_{h_l}(s)-Y^{\vv}_{h_l}(\eta_{h_{l}}(s)),\\
&~~~~~~~~~~~~~~~~~~f_{h_l}(X^{\vv}_{h_l}(\eta_{h_{l}}(s)), X^{\vv}_{h_l}(\eta_{h_{l}}(s-\tau)))\rangle\mbox{d}s\\
=:&pE_{31}(t)+pE_{32}(t),
\end{split}
\end{equation*}
where
\begin{equation*}
\begin{split}
Y^{\vv}_{h_l}(s)-Y^{\vv}_{h_l}(\eta_{h_{l}}(s))=&\int_{\eta_{h_{l}}(s)}^sf_{h_l}(X^{\vv}_{h_l}(\eta_{h_{l}}(u)), X^{\vv}_{h_l}(\eta_{h_{l}}(u-\tau)))\mbox{d}u\\
&+\int_{\eta_{h_{l}}(s)}^s\vv g(X^{\vv}_{h_l}(\eta_{h_{l}}(u)), X^{\vv}_{h_l}(\eta_{h_{l}}(u-\tau)))\mbox{d}W(u).
\end{split}
\end{equation*}
Due to \eqref{delbound} and the Young inequality,
\begin{equation*}
\begin{split}
&E_{31}(t)=\mathbb{E}\int_0^t[1+|Y^{\vv}_{h_l}(\eta_{h_{l}}(s))|^{2}]^{\frac{p-2}{2}}\bigg\langle\int_{\eta_{h_{l}}(s)}^sf_{h_l}(X^{\vv}_{h_l}(\eta_{h_{l}}(u)), X^{\vv}_{h_l}(\eta_{h_{l}}(u-\tau)))\mbox{d}u,\\
&~~~~~~~~~~~~~~~~~~~~~f_{h_l}(X^{\vv}_{h_l}(\eta_{h_{l}}(s)), X^{\vv}_{h_l}(\eta_{h_{l}}(s-\tau)))\bigg\rangle\mbox{d}s\\
&+\mathbb{E}\int_0^t[1+|Y^{\vv}_{h_l}(\eta_{h_{l}}(s))|^{2}]^{\frac{p-2}{2}}\bigg\langle\mathbb{E}\int_{\eta_{h_{l}}(s)}^s\vv g(X^{\vv}_{h_l}(\eta_{h_{l}}(u)),X^{\vv}_{h_l}(\eta_{h_{l}}(u-\tau)))\mbox{d}W(u)\bigg|_{\mathscr{F}_{\eta_{h_{l}}(s)}},\\
&~~~~~~~~~~~~~~~~~~~~~f_{h_l}(X^{\vv}_{h_l}(\eta_{h_{l}}(s)), X^{\vv}_{h_l}(\eta_{h_{l}}(s-\tau)))\bigg\rangle\mbox{d}s\\
&\le\mathbb{E}\int_0^t[1+|Y^{\vv}_{h_l}(\eta_{h_{l}}(s))|^{2}]^{\frac{p-2}{2}}\int_{\eta_{h_{l}}(s)}^s|f_{h_l}(X^{\vv}_{h_l}(\eta_{h_{l}}(u)), X^{\vv}_{h_l}(\eta_{h_{l}}(u-\tau)))|\mbox{d}u\\
&~~~~~~~~~~~~~~~~~~~~~|f_{h_l}(X^{\vv}_{h_l}(\eta_{h_{l}}(s)), X^{\vv}_{h_l}(\eta_{h_{l}}(s-\tau)))|\mbox{d}s\\
&\le h_{l-1}\mathbb{E}\int_0^t[1+|Y^{\vv}_{h_l}(\eta_{h_{l}}(s))|^{2}]^{\frac{p-2}{2}}|f_{h_l}(X^{\vv}_{h_l}(\eta_{h_{l}}(s)), X^{\vv}_{h_l}(\eta_{h_{l}}(s-\tau)))|^2\mbox{d}s\\
&\le Ch_{l-1}^{1-2\delta}\mathbb{E}\int_0^t(1+|X^{\vv}_{h_l}(\eta_{h_{l}}(s))|^p)\mbox{d}s+Ch_{l-1}^{1-2\delta}h_{l-1}^{(1-\delta)p}\\
&\le C+C\int_0^t\sup\limits_{0\le u\le s}\mathbb{E}|X^{\vv}_{h_l}(u)|^p\mbox{d}s.
\end{split}
\end{equation*}
Applying the It\^{o} formula again, we obtain
\begin{equation*}
\begin{split}
&[1+|Y^{\vv}_{h_l}(s)|^{2}]^{\frac{p-2}{2}}-[1+|Y^{\vv}_{h_l}(\eta_{h_{l}}(s))|^{2}]^{\frac{p-2}{2}}\le[1+|Y^{\vv}_{h_l}(0)|^{2}]^{\frac{p-2}{2}}-[1+|Y^{\vv}_{h_l}(\eta_{h_{l}}(0))|^{2}]^{\frac{p-2}{2}}\\
&~~~~+(p-2)\int_{\eta_{h_{l}}(s)}^s[1+|Y^{\vv}_{h_l}(u)|^{2}]^{\frac{p-4}{2}}\langle Y^{\vv}_{h_l}(u), f_{h_l}(X^{\vv}_{h_l}(\eta_{h_{l}}(u)), X^{\vv}_{h_l}(\eta_{h_{l}}(u-\tau)))\rangle\mbox{d}u\\
&~~~~+\frac{1}{2}(p-2)(p-3)\int_{\eta_{h_{l}}(s)}^s[1+|Y^{\vv}_{h_l}(u)|^{2}]^{\frac{p-4}{2}}|\vv g(X^{\vv}_{h_l}(\eta_{h_{l}}(u)), X^{\vv}_{h_l}(\eta_{h_{l}}(u-\tau)))|^2\mbox{d}u\\
&~~~~+(p-2)\int_{\eta_{h_{l}}(s)}^s[1+|Y^{\vv}_{h_l}(u)|^{2}]^{\frac{p-4}{2}}\langle Y^{\vv}_{h_l}(u), \vv g(X^{\vv}_{h_l}(\eta_{h_{l}}(u)), X^{\vv}_{h_l}(\eta_{h_{l}}(u-\tau)))\mbox{d}W(u)\rangle.
\end{split}
\end{equation*}
Hence,
\begin{equation*}
\begin{split}
E_{32}(t)\le&(p-2)\mathbb{E}\int_0^t\int_{\eta_{h_{l}}(s)}^s[1+|Y^{\vv}_{h_l}(u)|^{2}]^{\frac{p-4}{2}}\langle Y^{\vv}_{h_l}(u), f_{h_l}(X^{\vv}_{h_l}(\eta_{h_{l}}(u)), X^{\vv}_{h_l}(\eta_{h_{l}}(u-\tau)))\rangle\mbox{d}u\\
\quad&\times\langle Y^{\vv}_{h_l}(s)-Y^{\vv}_{h_l}(\eta_{h_{l}}(s)),f_{h_l}(X^{\vv}_{h_l}(\eta_{h_{l}}(s)), X^{\vv}_{h_l}(\eta_{h_{l}}(s-\tau)))\rangle\mbox{d}s\\
&+\frac{1}{2}(p-2)(p-3)\mathbb{E}\int_0^t\int_{\eta_{h_{l}}(s)}^s[1+|Y^{\vv}_{h_l}(u)|^{2}]^{\frac{p-4}{2}}|\vv g(X^{\vv}_{h_l}(\eta_{h_{l}}(u)), X^{\vv}_{h_l}(\eta_{h_{l}}(u-\tau)))|^2\mbox{d}u\\
\quad&\times\langle Y^{\vv}_{h_l}(s)-Y^{\vv}_{h_l}(\eta_{h_{l}}(s)),f_{h_l}(X^{\vv}_{h_l}(\eta_{h_{l}}(s)), X^{\vv}_{h_l}(\eta_{h_{l}}(s-\tau)))\rangle\mbox{d}s\\
&+(p-2)\mathbb{E}\int_0^t\int_{\eta_{h_{l}}(s)}^s[1+|Y^{\vv}_{h_l}(u)|^{2}]^{\frac{p-4}{2}}\langle Y^{\vv}_{h_l}(u), \vv g(X^{\vv}_{h_l}(\eta_{h_{l}}(u)), X^{\vv}_{h_l}(\eta_{h_{l}}(u-\tau)))\mbox{d}W(u)\rangle\\
\quad&\times\langle Y^{\vv}_{h_l}(s)-Y^{\vv}_{h_l}(\eta_{h_{l}}(s)),f_{h_l}(X^{\vv}_{h_l}(\eta_{h_{l}}(s)), X^{\vv}_{h_l}(\eta_{h_{l}}(s-\tau)))\rangle\mbox{d}s\\
=:&(p-2)E_{321}+\frac{1}{2}(p-2)(p-3)E_{322}+(p-2)E_{323}.
\end{split}
\end{equation*}
Using (H2), \eqref{delbound}, the Young inequality, the H\"{o}lder inequality and the Burkholder-Davis-Gundy (BDG) inequality, since $\delta\in(0,1/2]$, we compute
\begin{equation*}
\begin{split}
&E_{321}(t)\le\mathbb{E}\int_0^t\int_{\eta_{h_{l}}(s)}^s[1+|Y^{\vv}_{h_l}(u)|^{2}]^{\frac{p-4}{2}}\langle Y^{\vv}_{h_l}(u), f_{h_l}(X^{\vv}_{h_l}(\eta_{h_{l}}(u)), X^{\vv}_{h_l}(\eta_{h_{l}}(u-\tau)))\rangle\mbox{d}u\\
&\times\bigg\langle\int_{\eta_{h_{l}}(s)}^sf_{h_l}(X^{\vv}_{h_l}(\eta_{h_{l}}(u)), X^{\vv}_{h_l}(\eta_{h_{l}}(u-\tau)))\mbox{d}u,f_{h_l}(X^{\vv}_{h_l}(\eta_{h_{l}}(s)), X^{\vv}_{h_l}(\eta_{h_{l}}(s-\tau)))\bigg\rangle\mbox{d}s\\
&+\mathbb{E}\int_0^t\int_{\eta_{h_{l}}(s)}^s[1+|Y^{\vv}_{h_l}(u)|^{2}]^{\frac{p-4}{2}}\langle Y^{\vv}_{h_l}(u), f_{h_l}(X^{\vv}_{h_l}(\eta_{h_{l}}(u)), X^{\vv}_{h_l}(\eta_{h_{l}}(u-\tau)))\rangle\mbox{d}u\\
&\times\bigg\langle\int_{\eta_{h_{l}}(s)}^s\vv g(X^{\vv}_{h_l}(\eta_{h_{l}}(u)), X^{\vv}_{h_l}(\eta_{h_{l}}(u-\tau)))\mbox{d}W(u),f_{h_l}(X^{\vv}_{h_l}(\eta_{h_{l}}(s)), X^{\vv}_{h_l}(\eta_{h_{l}}(s-\tau)))\bigg\rangle\mbox{d}s\\
\le&h_{l-1}\mathbb{E}\int_0^t\int_{\eta_{h_{l}}(s)}^s[1+|Y^{\vv}_{h_l}(u)|^{2}]^{\frac{p-3}{2}}|f_{h_l}(X^{\vv}_{h_l}(\eta_{h_{l}}(s)), X^{\vv}_{h_l}(\eta_{h_{l}}(s-\tau)))|^3\mbox{d}u\mbox{d}s\\
&+C\mathbb{E}\int_0^t\bigg[\bigg(\int_{\eta_{h_{l}}(s)}^s[1+|Y^{\vv}_{h_l}(u)|^{2}]^{\frac{p-3}{2}}|f_{h_l}(X^{\vv}_{h_l}(\eta_{h_{l}}(u)), X^{\vv}_{h_l}(\eta_{h_{l}}(u-\tau)))|\mbox{d}u\\
&|f_{h_l}(X^{\vv}_{h_l}(\eta_{h_{l}}(s)), X^{\vv}_{h_l}(\eta_{h_{l}}(s-\tau)))|\bigg)^{\frac{p}{p-1}}+\left|\int_{\eta_{h_{l}}(s)}^s\vv g(X^{\vv}_{h_l}(\eta_{h_{l}}(u)), X^{\vv}_{h_l}(\eta_{h_{l}}(u-\tau)))\mbox{d}W(u)\right|^{p}\bigg]\mbox{d}s\\
\le&Ch_{l-1}^{2-3\delta}\mathbb{E}\int_0^t|X^{\vv}_{h_l}(s)|^{p}\mbox{d}s+Ch_{l-1}^{2-3\delta}h_{l-1}^{(1-\delta)p}\\
&+C\mathbb{E}\int_0^t\bigg(\int_{\eta_{h_{l}}(s)}^s[1+|Y^{\vv}_{h_l}(u)|^{2}]^{\frac{p-3}{2}}|f_{h_l}(X^{\vv}_{h_l}(\eta_{h_{l}}(s)), X^{\vv}_{h_l}(\eta_{h_{l}}(s-\tau)))|^2\mbox{d}u\bigg)^{\frac{p}{p-1}}\mbox{d}s\\
&+C\mathbb{E}\int_0^t\left(\int_{\eta_{h_{l}}(s)}^s|\vv g(X^{\vv}_{h_l}(\eta_{h_{l}}(u)), X^{\vv}_{h_l}(\eta_{h_{l}}(u-\tau)))|^2\mbox{d}u\right)^{\frac{p}{2}}\mbox{d}s\\
\le&C+C\int_0^t\sup\limits_{0\le u\le s}\mathbb{E}|X^{\vv}_{h_l}(u)|^{p}\mbox{d}s
+Ch_{l-1}^{(1-2\delta)\frac{p}{2}}\int_0^t\sup\limits_{0\le u\le s}\mathbb{E}|X^{\vv}_{h_l}(u)|^{p}\mbox{d}s\\
\le&C+C\int_0^t\sup\limits_{0\le u\le s}\mathbb{E}|X^{\vv}_{h_l}(u)|^{p}\mbox{d}s.
\end{split}
\end{equation*}
Using the same techniques in the way to the estimation of $E_{321}(t)$, we get
\begin{equation*}
E_{322}(t)\le C+C\int_0^t\sup\limits_{0\le u\le s}\mathbb{E}|X^{\vv}_{h_l}(u)|^{p}\mbox{d}s.
\end{equation*}
Furthermore, by (H2) and \eqref{delbound} again, we have
\begin{equation*}
\begin{split}
E_{323}(t)=&\mathbb{E}\int_0^t\int_{\eta_{h_{l}}(s)}^s[1+|Y^{\vv}_{h_l}(u)|^{2}]^{\frac{p-4}{2}}\langle Y^{\vv}_{h_l}(u),\vv g(X^{\vv}_{h_l}(\eta_{h_{l}}(u)), X^{\vv}_{h_l}(\eta_{h_{l}}(u-\tau)))\mbox{d}W(u)\rangle \\
&\times\bigg\langle\int_{\eta_{h_{l}}(s)}^sf_{h_l}(X^{\vv}_{h_l}(\eta_{h_{l}}(u)), X^{\vv}_{h_l}(\eta_{h_{l}}(u-\tau)))\mbox{d}u,
f_{h_l}(X^{\vv}_{h_l}(\eta_{h_{l}}(s)), X^{\vv}_{h_l}(\eta_{h_{l}}(s-\tau)))\bigg\rangle\mbox{d}s\\
&+\mathbb{E}\int_0^t\int_{\eta_{h_{l}}(s)}^s[1+|Y^{\vv}_{h_l}(u)|^{2}]^{\frac{p-4}{2}}\langle Y^{\vv}_{h_l}(u),\vv g(X^{\vv}_{h_l}(\eta_{h_{l}}(u)), X^{\vv}_{h_l}(\eta_{h_{l}}(u-\tau)))\mbox{d}W(u)\rangle \\
&\times\bigg\langle\int_{\eta_{h_{l}}(s)}^s\vv g(X^{\vv}_{h_l}(\eta_{h_{l}}(u)), X^{\vv}_{h_l}(\eta_{h_{l}}(u-\tau)))\mbox{d}W(u),
f_{h_l}(X^{\vv}_{h_l}(\eta_{h_{l}}(s)), X^{\vv}_{h_l}(\eta_{h_{l}}(s-\tau)))\bigg\rangle\mbox{d}s\\
=&\mathbb{E}\int_0^t\int_{\eta_{h_{l}}(s)}^s[1+|Y^{\vv}_{h_l}(u)|^{2}]^{\frac{p-4}{2}}\langle Y^{\vv}_{h_l}(u),\vv g(X^{\vv}_{h_l}(\eta_{h_{l}}(u)), X^{\vv}_{h_l}(\eta_{h_{l}}(u-\tau)))\mbox{d}W(u)\rangle \\
&\times\bigg\langle\int_{\eta_{h_{l}}(s)}^s\vv g(X^{\vv}_{h_l}(\eta_{h_{l}}(u)), X^{\vv}_{h_l}(\eta_{h_{l}}(u-\tau)))\mbox{d}W(u),
f_{h_l}(X^{\vv}_{h_l}(\eta_{h_{l}}(s)), X^{\vv}_{h_l}(\eta_{h_{l}}(s-\tau)))\bigg\rangle\mbox{d}s\\
\le&\mathbb{E}\int_0^t\int_{\eta_{h_{l}}(s)}^s[1+|Y^{\vv}_{h_l}(u)|^{2}]^{\frac{p-3}{2}}|\vv g(X^{\vv}_{h_l}(\eta_{h_{l}}(u)), X^{\vv}_{h_l}(\eta_{h_{l}}(u-\tau)))|^2\mbox{d}u\\
&~~~~~~~~~~~|f_{h_l}(X^{\vv}_{h_l}(\eta_{h_{l}}(s)), X^{\vv}_{h_l}(\eta_{h_{l}}(s-\tau)))|\mbox{d}s\\
\le&Ch_{l-1}^{1-\delta}\mathbb{E}\int_0^t(|X^{\vv}_{h_l}(s)|^{p}+|X^{\vv}_{h_l}(\eta_{h_{l}}(s))|^{p}+|X^{\vv}_{h_l}(\eta_{h_{l}}(s-\tau)))|^{p})\mbox{d}s+Ch_{l-1}^{1-\delta}h_{l-1}^{(1-\delta)p}\\
\le& C+C\int_0^t\sup\limits_{0\le u\le s}\mathbb{E}|X^{\vv}_{h_l}(u)|^{p}\mbox{d}s.
\end{split}
\end{equation*}
By sorting these equations, we conclude that
\begin{equation*}\label{e3t}
E_3(t)\le C+C\int_0^t\sup\limits_{0\le u\le s}\mathbb{E}|X^{\vv}_{h_l}(u)|^{p}\mbox{d}s.
\end{equation*}
Thus, the estimation of $E_1(t)-E_3(t)$ results in
\begin{equation}\label{tildeyt}
\begin{split}
\sup\limits_{0\le u\le t}\mathbb{E}|Y^{\vv}_{h_l}(u)|^p\le\sup\limits_{0\le u\le t}\mathbb{E}[1+|Y^{\vv}_{h_l}(u)|^{2}]^{\frac{p}{2}}\le C+C\int_0^t\sup\limits_{0\le u\le s}\mathbb{E}|X^{\vv}_{h_l}(u)|^p\mbox{d}s.
\end{split}
\end{equation}
By the relationship between $X^{\vv}_{h_l}(t)$ and $Y^{\vv}_{h_l}(t)$, it is easy to derive from \eqref{delbound} that
\begin{equation*}
\begin{split}
\sup\limits_{0\le u\le t}\mathbb{E}|X^{\vv}_{h_l}(u)|^p\le Ch_{l-1}^{(1-\delta)p}+\sup\limits_{0\le u\le t}\mathbb{E}|Y^{\vv}_{h_l}(u)|^p
&\le C+C\int_0^t\sup\limits_{0\le u\le s}\mathbb{E}|X^{\vv}_{h_l}(u)|^p\mbox{d}s.
\end{split}
\end{equation*}
Finally, the desired result follows by the Gronwall inequality. $\Box$

Let $Z_{h_l}$ be the deterministic solution to
\begin{equation*}
Z_{h_l}(t)-\theta f_{h_l}(Z_{h_l}(t), Z_{h_l}(t-\tau)){h_l}=\xi(0)-\theta f_{h_l}(\xi(0),\xi(-\tau)){h_l}+\int_0^tf_{h_l}(Z_{h_l}(\eta_{h_l}(s)), Z_{h_l}(\eta_{h_l}(s-\tau))){\mbox d}s,
\end{equation*}
which is the corresponding theta EM approximation to the ordinary differential delay equation obtained from \eqref{2.0} when $\vv=0$.

\begin{lemma}\label{1sun}
Let assumptions (H1) and (H2) hold. Then, for any $T>0$ and $p\ge 2,$ we have
$$
\sup_{0\le t\le T}\E|Z_{h_{l}}(t)|^p\le C,
$$
and
$$
\sup_{0\le n<M^{l-1}, 1\le k\le M}\E |Z_{h_{l}}(t_n^k)-Z_{h_l}(t_n)|^p \le CM^{(1-\delta)p}h_l^{(1-\delta)p}.
$$
\end{lemma}
{\bf Proof.} Following the proof of Lemma \ref{1pmoment}, the first part is obvious. Denote by $\bar{Z}_{h_l}(t)=Z_{h_l}(t)-\theta f_{h_l}(Z_{h_l}(t), Z_{h_l}(t-\tau))h_l$. For any $n\in \{0, 1, \ldots, M^{l-1}-1\}$ and $k\in \{1, \ldots, M\}$, by \eqref{delbound}, we have
\begin{equation*}
\begin{split}
\E|\bar{Z}_{h_l}(t_n^k)-\bar{Z}_{h_l}(t_n)|^p\le |kh_l|^{p-1}\E\int_{t_n}^{t_n^k}|f_{h_l}(Z_{h_l}(\eta_{h_l}(s)), Z_{h_l}(\eta_{h_l}(s-\tau)))|^p{\mbox d}s\le CM^{(1-\delta)p}h_l^{(1-\delta)p}.
\end{split}
\end{equation*}
On the other side, we see
\begin{equation*}
\begin{split}
\E|Z_{h_l}(t_n^k)-Z_{h_l}(t_n)|^p\le C\E|\bar{Z}_{h_l}(t_n^k)-\bar{Z}_{h_l}(t_n)|^p+CM^{(1-\delta)p}h_l^{(1-\delta)p}.
\end{split}
\end{equation*}
Thus, the desired assertion follows. \hfill $\Box$

\begin{lemma}\label{1miss}
Let assumptions (H1) and (H2) hold. Then, for any $p>0,$ we have
\begin{equation*}
\sup_{0\le n< M^{l-1}, 1\le k\le M}\E [|X^\vv_{h_{l}}(t_n^k)-X_{h_l}^\vv(t_n)|^p]\le CM^{(1-\delta)p} h_l^{(1-\delta)p}+C\vv^pM^{p/2}h_l^{p/2}.
\end{equation*}
\end{lemma}
{\bf Proof.}
For $n\in \{0, 1, \ldots, M^{l-1}-1\}$ and $k\in \{1, \ldots, M\}$, we see
\begin{equation*}
\begin{split}
&Y^{\vv}_{h_l}(t_n^k)=Y^{\vv}_{h_l}(t_{n})+\int_{t_{n}}^{t_n^k}f_{h_l}(X^{\vv}_{h_l}(\eta_{h_{l}}(s)), X^{\vv}_{h_l}(\eta_{h_{l}}(s-\tau))){\mbox d}s\\
&+\vv\int_{t_{n}}^{t_n^k}g(X^{\vv}_{h_l}(\eta_{h_{l}}(s)), X^{\vv}_{h_l}(\eta_{h_{l}}(s-\tau))){\mbox d}W(s).
\end{split}
\end{equation*}
By the elementary inequality $|a+b|^p\le 2^{p-1}(|a|^p+|b|^p), p\ge 1$, we compute for $p\ge1$
\begin{equation*}
\begin{split}
\mathbb{E}|Y^{\vv}_{h_l}(t_n^k)-Y^{\vv}_{h_l}(t_{n})|^p&\le 2^{p-1}\mathbb{E}\left|\int_{t_n}^{t_n^k} f_{h_l}(X^{\vv}_{h_l}(\eta_{h_{l}}(s)), X^{\vv}_{h_l}(\eta_{h_{l}}(s-\tau)))\mbox{d}s\right|^p\\
&+2^{p-1}\vv^p\mathbb{E}\left|\int_{t_n}^{t_n^k} g(X^{\vv}_{h_l}(\eta_{h_{l}}(s)), X^{\vv}_{h_l}(\eta_{h_{l}}(s-\tau))){\mbox d}W(s)\right|^p.
\end{split}
\end{equation*}
With (H2), \eqref{delbound}, Lemma \ref{1pmoment}, the H\"{o}lder inequality and the BDG inequality, we derive
\begin{equation}\label{1niu}
\begin{split}
&\mathbb{E}|Y^{\vv}_{h_l}(t_n^k)-Y^{\vv}_{h_l}(t_{n})|^p
\le 2^{p-1}M^{p-1}h_{l}^{p-1}\mathbb{E}\int_{t_n}^{t_n^k} \left|f_{h_l}(X^{\vv}_{h_l}(\eta_{h_{l}}(s)), X^{\vv}_{h_l}(\eta_{h_{l}}(s-\tau)))\right|^p\mbox{d}s\\
&+C\vv^pM^{\frac{p}{2}-1}h_{l}^{\frac{p}{2}-1}\mathbb{E}\int_{t_n}^{t_n^k} \left|g(X^{\vv}_{h_l}(\eta_{h_{l}}(s)), X^{\vv}_{h_l}(\eta_{h_{l}}(s-\tau)))\right|^p\mbox{d}s\\
\le&CM^{(1-\delta)p}h_{l}^{(1-\delta)p}+C\vv^pM^{p/2}h_{l}^{p/2}.
\end{split}
\end{equation}
Since we have
\begin{equation*}
\begin{split}
X^{\vv}_{h_l}(t_n^k)-X^{\vv}_{h_l}(t_{n})=&Y^{\vv}_{h_l}(t_n^k)-Y^{\vv}_{h_l}(t_{n})+\theta f_{h_l}(X^{\vv}_{h_{l}}(t_n^k),X^{\vv}_{h_{l}}(t_{n}^k-m_lh_l))h_l\\
&-\theta f_{h_l}(X^{\vv}_{h_{l}}(t_{n}),X^{\vv}_{h_{l}}(t_{n}-m_lh_l))h_l.
\end{split}
\end{equation*}
This combines with \eqref{1niu} lead to
\begin{equation*}
\begin{split}
\mathbb{E}|X^{\vv}_{h_l}(t_n^k)-X^{\vv}_{h_l}(t_{n})|^p\le&C\mathbb{E}|Y^{\vv}_{h_l}(t_n^k)-Y^{\vv}_{h_l}(t_{n})|^p+CM^{(1-\delta)p}h_l^{(1-\delta)p}\\
\le&CM^{(1-\delta)p}h_{l}^{(1-\delta)p}+C\vv^pM^{p/2}h_{l}^{p/2}.
\end{split}
\end{equation*}
The desired result then follows for $p\ge1$. Finally, one can use the Young inequality to get the results for $p\in(0,1)$.\hfill $\Box$

\begin{lemma}\label{1gaosi}
Let assumptions (H1) and (H2) hold. Then, for any $T>0$ and $p>0,$ we have
$$
\E \left[\sup_{0\le n<M^l}\sup_{nh_l\le t<(n+1)h_l}|X_{h_{l}}^\vv(t)-X_{h_{l}}^\vv(nh_l)|^p\right] \le  Ch_l^{(1-\delta)p}+C\vv^ph_l^{p/2},
$$
and
$$
\E \left[\sup_{0\le n<M^l}\sup_{nh_l\le t<(n+1)h_l}|Z_{h_{l}}(t)-Z_{h_{l}}(nh_l)|^p\right] \le Ch_l^{(1-\delta)p}.
$$
\end{lemma}
{\bf Proof.} Similar to the proof of Lemma \ref{1miss}, the result is obvious. \hfill $\Box$

\begin{lemma}\label{1tech}
Let assumptions (H1) and (H2) hold. Then, for any $T>0$ and $p\ge 2,$ we have
$$
\E \left[\sup_{0\le t\le T}|X_{h_{l}}^\vv(t)-Z_{h_l}(t)|^p\right] \le CM^{\delta p/2}h_l^{\delta p/2}+C\vv^{\frac{p}{2}}M^{{\frac{1-2\delta}{4}p}}h_l^{{\frac{1-2\delta}{4}p}}+C\vv^p,
$$
and
$$
\E \left[\sup_{0\le t\le T}|X_{h_{l-1}}^\vv(t)-Z_{h_{l-1}}(t)|^p\right]\le CM^{\delta p/2}h_{l-1}^{\delta p/2}+C\vv^{\frac{p}{2}}M^{{\frac{1-2\delta}{4}p}}h_{l-1}^{{\frac{1-2\delta}{4}p}}+C\vv^p.
$$
\end{lemma}
{\bf Proof.} By the definition of $Y^{\vv}_{h_l}(t)$ and $\bar{Z}_{h_l}(t)$,
\begin{equation*}
\begin{split}
Y^{\vv}_{h_l}(t)-\bar{Z}_{h_l}(t)=&\int_0^t[f_{h_l}(X^{\vv}_{h_l}(\eta_{h_{l}}(s)), X^{\vv}_{h_l}(\eta_{h_{l}}(s-\tau)))-f_{h_l}(Z_{h_l}(\eta_{h_l}(s)), Z_{h_l}(\eta_{h_l}(s-\tau)))]{\mbox d}s\\
&+\vv\int_0^tg(X^{\vv}_{h_l}(\eta_{h_{l}}(s)), X^{\vv}_{h_l}(\eta_{h_{l}}(s-\tau))){\mbox d}W(s),
\end{split}
\end{equation*}
thus, by the It\^{o} formula,
\begin{equation*}
\begin{split}
&|Y^{\vv}_{h_l}(t)-\bar{Z}_{h_l}(t)|^p\le p\int_0^t|Y^{\vv}_{h_l}(s)-\bar{Z}_{h_l}(s)|^{p-2}\big\langle Y^{\vv}_{h_l}(s)-\bar{Z}_{h_l}(s), \\
&f_{h_l}(X^{\vv}_{h_l}(\eta_{h_{l}}(s)), X^{\vv}_{h_l}(\eta_{h_{l}}(s-\tau)))-f_{h_l}(Z_{h_l}(\eta_{h_l}(s)), Z_{h_l}(\eta_{h_l}(s-\tau)))\big\rangle{\mbox d}s\\
&+\frac{p(p-1)}{2}\int_0^t|Y^{\vv}_{h_l}(s)-\bar{Z}_{h_l}(s)|^{p-2}|\vv g(X^{\vv}_{h_l}(\eta_{h_{l}}(s)), X^{\vv}_{h_l}(\eta_{h_{l}}(s-\tau)))|^2{\mbox d}s\\
&+p\int_0^t|Y^{\vv}_{h_l}(s)-\bar{Z}_{h_l}(s)|^{p-2}\big\langle Y^{\vv}_{h_l}(s)-\bar{Z}_{h_l}(s),\vv g(X^{\vv}_{h_l}(\eta_{h_{l}}(s)), X^{\vv}_{h_l}(\eta_{h_{l}}(s-\tau))){\mbox d}W(s)\big\rangle\\
\le&H_1(t)+H_2(t)+H_3(t)+H_4(t)+H_5(t)+H_6(t)+H_7(t)+H_8(t),
\end{split}
\end{equation*}
where
%\begin{enumerate}
%\item[]
\begin{equation*}
\begin{split}
H_1(t)=&p\int_0^t|Y^{\vv}_{h_l}(s)-\bar{Z}_{h_l}(s)|^{p-2}\big\langle X^{\vv}_{h_l}(s)-X^{\vv}_{h_l}(\eta_{h_l}(s)), \\
&f_{h_l}(X^{\vv}_{h_l}(\eta_{h_{l}}(s)), X^{\vv}_{h_l}(\eta_{h_{l}}(s-\tau)))-f_{h_l}(Z_{h_l}(\eta_{h_l}(s)), Z_{h_l}(\eta_{h_l}(s-\tau)))\big\rangle{\mbox d}s,
\end{split}
\end{equation*}
%\item[]
\begin{equation*}
\begin{split}
H_2(t)=&p\int_0^t|Y^{\vv}_{h_l}(s)-\bar{Z}_{h_l}(s)|^{p-2}\big\langle X^{\vv}_{h_l}(\eta_{h_l}(s))-Z_{h_l}(\eta_{h_l}(s)), \\
&f_{h_l}(X^{\vv}_{h_l}(\eta_{h_{l}}(s)), X^{\vv}_{h_l}(\eta_{h_{l}}(s-\tau)))-f(X^{\vv}_{h_l}(\eta_{h_{l}}(s)), X^{\vv}_{h_l}(\eta_{h_{l}}(s-\tau)))\big\rangle{\mbox d}s,
\end{split}
\end{equation*}
%\item[]
\begin{equation*}
\begin{split}
H_3(t)=&p\int_0^t|Y^{\vv}_{h_l}(s)-\bar{Z}_{h_l}(s)|^{p-2}\big\langle X^{\vv}_{h_l}(\eta_{h_l}(s))-Z_{h_l}(\eta_{h_l}(s)), \\
&f(X^{\vv}_{h_l}(\eta_{h_{l}}(s)), X^{\vv}_{h_l}(\eta_{h_{l}}(s-\tau)))-f(Z_{h_l}(\eta_{h_l}(s)), Z_{h_l}(\eta_{h_l}(s-\tau)))\big\rangle{\mbox d}s,
\end{split}
\end{equation*}
%\item[]
\begin{equation*}
\begin{split}
H_4(t)=&p\int_0^t|Y^{\vv}_{h_l}(s)-\bar{Z}_{h_l}(s)|^{p-2}\big\langle X^{\vv}_{h_l}(\eta_{h_l}(s))-Z_{h_l}(\eta_{h_l}(s)), \\
&f(Z_{h_l}(\eta_{h_l}(s)), Z_{h_l}(\eta_{h_l}(s-\tau)))-f_{h_l}(Z_{h_l}(\eta_{h_l}(s)), Z_{h_l}(\eta_{h_l}(s-\tau)))\big\rangle{\mbox d}s,
\end{split}
\end{equation*}
%\item[]
\begin{equation*}
\begin{split}
H_5(t)=&p\int_0^t|Y^{\vv}_{h_l}(s)-\bar{Z}_{h_l}(s)|^{p-2}\big\langle Z_{h_l}(\eta_{h_l}(s))-Z_{h_l}(s),\\
&f_{h_l}(X^{\vv}_{h_l}(\eta_{h_{l}}(s)), X^{\vv}_{h_l}(\eta_{h_{l}}(s-\tau)))-f_{h_l}(Z_{h_l}(\eta_{h_l}(s)), Z_{h_l}(\eta_{h_l}(s-\tau)))\big\rangle{\mbox d}s,
\end{split}
\end{equation*}
%\item[]
\begin{equation*}
\begin{split}
H_6(t)=&-\theta ph_l\int_0^t|Y^{\vv}_{h_l}(s)-\bar{Z}_{h_l}(s)|^{p-2}\big\langle f_{h_l}(X^{\vv}_{h_l}(s), X^{\vv}_{h_l}(s-\tau))-f_{h_l}(Z_{h_l}(s), Z_{h_l}(s-\tau)),\\
&f_{h_l}(X^{\vv}_{h_l}(\eta_{h_{l}}(s)), X^{\vv}_{h_l}(\eta_{h_{l}}(s-\tau)))-f_{h_l}(Z_{h_l}(\eta_{h_l}(s)), Z_{h_l}(\eta_{h_l}(s-\tau)))\big\rangle{\mbox d}s,
\end{split}
\end{equation*}
%\item[]
\begin{equation*}
\begin{split}
H_7(t)=&\frac{p(p-1)}{2}\int_0^t|Y^{\vv}_{h_l}(s)-\bar{Z}_{h_l}(s)|^{p-2}|\vv g(X^{\vv}_{h_l}(\eta_{h_{l}}(s)), X^{\vv}_{h_l}(\eta_{h_{l}}(s-\tau)))|^2{\mbox d}s,
\end{split}
\end{equation*}
%\item[]
\begin{equation*}
\begin{split}
H_8(t)=&p\int_0^t|Y^{\vv}_{h_l}(s)-\bar{Z}_{h_l}(s)|^{p-2}\big\langle Y^{\vv}_{h_l}(s)-\bar{Z}_{h_l}(s),\vv g(X^{\vv}_{h_l}(\eta_{h_{l}}(s)), X^{\vv}_{h_l}(\eta_{h_{l}}(s-\tau))){\mbox d}W(s)\big\rangle.
\end{split}
\end{equation*}
%\end{enumerate}
By the H\"{o}lder inequality, \eqref{delbound} and Lemma \ref{1gaosi}, we get
\begin{equation*}
\begin{split}
&\E\left(\sup\limits_{0\le u\le t}|H_1(u)|^p\right)\le C\E\int_0^t|Y^{\vv}_{h_l}(s)-\bar{Z}_{h_l}(s)|^{p}{\mbox d}s+C\E\int_0^t|X^{\vv}_{h_l}(s)-X^{\vv}_{h_l}(\eta_{h_l}(s))|^{p/2}\\
&|f_{h_l}(X^{\vv}_{h_l}(\eta_{h_{l}}(s)), X^{\vv}_{h_l}(\eta_{h_{l}}(s-\tau)))-f_{h_l}(Z_{h_l}(\eta_{h_l}(s)), Z_{h_l}(\eta_{h_l}(s-\tau)))|^{p/2}{\mbox d}s\\
&\le C\int_0^t\E \left[\sup_{0\le u\le s}|X_{h_{l}}^\vv(u)-Z_{h_l}(u)|^{p}\right]{\mbox d}s+CM^{\frac{1-2\delta}{2}p}h_l^{\frac{1-2\delta}{2}p}+C\vv^{\frac{p}{2}}M^{\frac{1-2\delta}{4}p}h_l^{\frac{1-2\delta}{4}p}.
\end{split}
\end{equation*}
By \eqref{delbound}, \eqref{del}, Lemma \ref{1pmoment}, Lemma \ref{1sun} and the H\"{o}lder inequality again,
\begin{equation*}
\begin{split}
&\E\left(\sup\limits_{0\le u\le t}|H_2(u)|^p\right)\le C\E\int_0^t|Y^{\vv}_{h_l}(s)-\bar{Z}_{h_l}(s)|^{p}{\mbox d}s+C\E\int_0^t|X^{\vv}_{h_l}(\eta_{h_l}(s))-Z_{h_l}(\eta_{h_l}(s))|^{p/2}\\
&|f_{h_l}(X^{\vv}_{h_l}(\eta_{h_{l}}(s)), X^{\vv}_{h_l}(\eta_{h_{l}}(s-\tau)))-f(X^{\vv}_{h_l}(\eta_{h_{l}}(s)), X^{\vv}_{h_l}(\eta_{h_{l}}(s-\tau))|^{p/2}{\mbox d}s\\
&\le C\int_0^t\E \left[\sup_{0\le u\le s}|X_{h_{l}}^\vv(u)-Z_{h_l}(u)|^{p}\right]{\mbox d}s+CM^{(1-\delta)p}h_l^{(1-\delta)p}+CM^{\delta p/2}h_l^{\delta p/2}.
\end{split}
\end{equation*}
By assumption (H1) and the H\"{o}lder inequality,
\begin{equation*}
\begin{split}
&\E\left(\sup\limits_{0\le u\le t}|H_3(u)|^p\right)\le C\E\int_0^t|Y^{\vv}_{h_l}(s)-\bar{Z}_{h_l}(s)|^{p}{\mbox d}s\\
&+C\E\int_0^t|X^{\vv}_{h_l}(\eta_{h_l}(s))-Z_{h_l}(\eta_{h_l}(s))|^{p}{\mbox d}s+C\E\int_0^t|X^{\vv}_{h_l}(\eta_{h_l}(s-\tau))-Z_{h_l}(\eta_{h_l}(s-\tau))|^{p}{\mbox d}s\\
\le &C\int_0^t\E \left[\sup_{0\le u\le s}|X_{h_{l}}^\vv(u)-Z_{h_l}(u)|^{p}\right]{\mbox d}s+CM^{(1-\delta)p}h_l^{(1-\delta)p}.
\end{split}
\end{equation*}
Similar to the estimation of $H_2(t)$, we derive
\begin{equation*}
\begin{split}
&\E\left(\sup\limits_{0\le u\le t}|H_4(u)|^p\right)\le C\int_0^t\E \left[\sup_{0\le u\le s}|X_{h_{l}}^\vv(u)-Z_{h_l}(u)|^p\right]{\mbox d}s+CM^{(1-\delta)p}h_l^{(1-\delta)p}+CM^{\delta p/2}h_l^{\delta p/2}.
\end{split}
\end{equation*}
With \eqref{delbound} and Lemma \ref{1gaosi}, we arrive at
\begin{equation*}
\begin{split}
&\E\left(\sup\limits_{0\le u\le t}|H_5(u)|^p\right)\le C\int_0^t\E \left[\sup_{0\le u\le s}|X_{h_{l}}^\vv(u)-Z_{h_l}(u)|^p\right]{\mbox d}s+CM^{\frac{1-2\delta}{2}p}h_l^{\frac{1-2\delta}{2}p}.
\end{split}
\end{equation*}
By \eqref{delbound} again, it is easy to see
\begin{equation*}
\begin{split}
&\E\left(\sup\limits_{0\le u\le t}|H_6(u)|^p\right)\le C\int_0^t\E \left[\sup_{0\le u\le s}|X_{h_{l}}^\vv(u)-Z_{h_l}(u)|^p\right]{\mbox d}s+CM^{(1-\delta)p}h_l^{(1-\delta)p}.
\end{split}
\end{equation*}
Then, (H2), \eqref{delbound}, Lemma \ref{1pmoment} and the H\"{o}lder inequality lead to
\begin{equation*}
\begin{split}
\E\left(\sup\limits_{0\le u\le t}|H_7(u)|^p\right)\le &C\E\int_0^t|Y^{\vv}_{h_l}(s)-\bar{Z}_{h_l}(s)|^{p}{\mbox d}s+C\vv^p\\
&\le C\int_0^t\E \left[\sup_{0\le u\le s}|X_{h_{l}}^\vv(u)-Z_{h_l}(u)|^p\right]{\mbox d}s+CM^{(1-\delta)p}h_l^{(1-\delta)p}+C\vv^p.
\end{split}
\end{equation*}
Moreover, with (H2), \eqref{delbound}, Lemma \ref{1pmoment} and the BDG inequality,
\begin{equation*}
\begin{split}
\E\left(\sup\limits_{0\le u\le t}|H_8(u)|^p\right)\le &\frac{1}{4}\E \left[\sup_{0\le s\le t}|Y_{h_{l}}^\vv(s)-\bar{Z}_{h_l}(s)|^p\right]+CM^{(1-\delta)p}h_l^{(1-\delta)p}+C\vv^p.
\end{split}
\end{equation*}
Sorting the above inequations together leads to
\begin{equation*}
\begin{split}
&\E \left[\sup_{0\le u\le t}|Y_{h_{l}}^\vv(u)-\bar{Z}_{h_l}(u)|^p\right]\le C\int_0^t\E \left[\sup_{0\le u\le s}|X_{h_{l}}^\vv(u)-Z_{h_l}(u)|^p\right]{\mbox d}s\\
&+ C\vv^p+C\vv^{\frac{p}{2}}M^{{\frac{1-2\delta}{4}p}}h_l^{{\frac{1-2\delta}{4}p}}+CM^{\delta p/2}h_l^{\delta p/2}.
\end{split}
\end{equation*}
Since we have
\begin{equation*}
\begin{split}
&\E \left[\sup_{0\le u\le t}|X_{h_{l}}^\vv(u)-Z_{h_l}(u)|^p\right]\le \E \left[\sup_{0\le u\le t}|Y_{h_{l}}^\vv(u)-\bar{Z}_{h_l}(u)|^p\right]+CM^{(1-\delta)p}h_l^{(1-\delta)p}.
\end{split}
\end{equation*}
The Gronwall inequality then leads to the desired result. By the same technique, the second part can be verified. \hfill $\Box$

\begin{theorem}\label{1election}
Let assumptions (H1) and (H2) hold. Then
\begin{equation*}
\sup_{0\le t\le T}\E|X^\vv_{h_{l}}(t)-X_{h_{l-1}}^\vv(t)|^2\le Ch_{l-1}^{1-2\delta}+Ch_{l-2}^{2\delta}+C\vv^2h_{l-1}.
\end{equation*}
\end{theorem}
{\bf Proof.} For any $t\in[0, T]$ , by \eqref{1c1} and \eqref{1c2}, we get
\begin{equation*}
\begin{split}
&Y^\vv_{h_{l}}(t)-Y_{h_{l-1}}^\vv(t)=Y^\vv_{h_{l}}(0)-Y_{h_{l-1}}^\vv(0)\\
&+\int_0^{t}[f_{h_l}(X^{\vv}_{h_l}(\eta_{h_l}(s)), X^{\vv}_{h_l}(\eta_{h_l}(s-\tau)))-f_{h_{l-1}}(X^{\vv}_{h_{l-1}}(\eta_{h_{l-1}}(s)), X^{\vv}_{h_{l-1}}(\eta_{h_{l-1}}(s-\tau)))]{\mbox d}s\\
&+\vv\int_0^{t}[g(X^{\vv}_{h_l}(\eta_{h_l}(s)), X^{\vv}_{h_l}(\eta_{h_l}(s-\tau)))-g(X^{\vv}_{h_{l-1}}(\eta_{h_{l-1}}(s)), X^{\vv}_{h_{l-1}}(\eta_{h_{l-1}}(s-\tau)))]{\mbox d}W(s).
\end{split}
\end{equation*}
By the It\^{o} formula, we get
\begin{equation*}
\begin{split}
&\E|Y^\vv_{h_{l}}(t)-Y_{h_{l-1}}^\vv(t)|^2\le 2\E|Y^\vv_{h_{l}}(0)-Y_{h_{l-1}}^\vv(0)|^2+2\E\int_0^{t}\langle Y^\vv_{h_{l}}(s)-Y_{h_{l-1}}^\vv(s),\\
& f_{h_l}(X^{\vv}_{h_l}(\eta_{h_l}(s)), X^{\vv}_{h_l}(\eta_{h_l}(s-\tau)))-f_{h_{l-1}}(X^{\vv}_{h_{l-1}}(\eta_{h_{l-1}}(s)), X^{\vv}_{h_{l-1}}(\eta_{h_{l-1}}(s-\tau)))\rangle{\mbox d}s\\
&+C\vv^2\E\int_0^{t}|g(X^{\vv}_{h_l}(\eta_{h_l}(s)), X^{\vv}_{h_l}(\eta_{h_l}(s-\tau)))-g(X^{\vv}_{h_{l-1}}(\eta_{h_{l-1}}(s)), X^{\vv}_{h_{l-1}}(\eta_{h_{l-1}}(s-\tau)))|^2{\mbox d}s\\
= &2\E|Y^\vv_{h_{l}}(0)-Y_{h_{l-1}}^\vv(0)|^2+2\E\int_0^{t}\langle Y^\vv_{h_{l}}(s)-Y_{h_{l-1}}^\vv(s),\\
&~~~~~f_{h_l}(X^{\vv}_{h_l}(\eta_{h_l}(s)), X^{\vv}_{h_l}(\eta_{h_l}(s-\tau)))-f(X^{\vv}_{h_l}(\eta_{h_l}(s)), X^{\vv}_{h_l}(\eta_{h_l}(s-\tau)))\rangle{\mbox d}s\\
&+2\E\int_0^{t}\langle Y^\vv_{h_{l}}(s)-Y_{h_{l-1}}^\vv(s), f(X^{\vv}_{h_l}(\eta_{h_l}(s)), X^{\vv}_{h_l}(\eta_{h_l}(s-\tau)))\\
&~~~-f_{h_{l-1}}(X^{\vv}_{h_l}(\eta_{h_l}(s)), X^{\vv}_{h_l}(\eta_{h_l}(s-\tau)))\rangle{\mbox d}s+2\E\int_0^{t}\langle X^{\vv}_{h_l}(\eta_{h_l}(s))-X^{\vv}_{h_{l-1}}(\eta_{h_{l-1}}(s)),\\
&~~~~~f_{h_{l-1}}(X^{\vv}_{h_l}(\eta_{h_l}(s)), X^{\vv}_{h_l}(\eta_{h_l}(s-\tau)))-f_{h_{l-1}}(X^{\vv}_{h_{l-1}}(\eta_{h_{l-1}}(s)), X^{\vv}_{h_{l-1}}(\eta_{h_{l-1}}(s-\tau)))\rangle{\mbox d}s\\
&+2\E\int_0^{t}\langle X^\vv_{h_{l}}(s)-X^{\vv}_{h_l}(\eta_{h_l}(s))-X_{h_{l-1}}^\vv(s)+X^{\vv}_{h_{l-1}}(\eta_{h_{l-1}}(s)),\\
&~~~~~f_{h_{l-1}}(X^{\vv}_{h_l}(\eta_{h_l}(s)), X^{\vv}_{h_l}(\eta_{h_l}(s-\tau)))-f_{h_{l-1}}(X^{\vv}_{h_{l-1}}(\eta_{h_{l-1}}(s)), X^{\vv}_{h_{l-1}}(\eta_{h_{l-1}}(s-\tau)))\rangle{\mbox d}s\\
&-2\theta h_l\E\int_0^{t}\langle f_{h_l}(X^{\vv}_{h_l}(s), X^{\vv}_{h_l}(s-\tau))-f_{h_{l-1}}(X^{\vv}_{h_{l-1}}(s), X^{\vv}_{h_{l-1}}(s-\tau)),\\
&~~~~~f_{h_l}(X^{\vv}_{h_l}(\eta_{h_l}(s)), X^{\vv}_{h_l}(\eta_{h_l}(s-\tau)))-f_{h_{l-1}}(X^{\vv}_{h_{l-1}}(\eta_{h_{l-1}}(s)), X^{\vv}_{h_{l-1}}(\eta_{h_{l-1}}(s-\tau)))\rangle{\mbox d}s\\
&+C\vv^2\E\int_0^{t}|g(X^{\vv}_{h_l}(\eta_{h_l}(s)), X^{\vv}_{h_l}(\eta_{h_l}(s-\tau)))-g(X^{\vv}_{h_{l-1}}(\eta_{h_{l-1}}(s)), X^{\vv}_{h_{l-1}}(\eta_{h_{l-1}}(s-\tau)))|^2{\mbox d}s\\
&\le Ch_{l-1}^{2-2\delta}+J_1(t)+J_2(t)+J_3(t)+J_4(t)+J_5(t)+J_6(t).
\end{split}
\end{equation*}
By \eqref{delbound}, \eqref{del} and Lemma \ref{1pmoment} we see
\begin{equation*}
\begin{split}
&J_1(t)\le C\E\int_0^{t}|X^\vv_{h_{l}}(s)-X_{h_{l-1}}^\vv(s)|^2{\mbox d}s+Ch_{l}^{2-2\delta}+Ch_{l-1}^{2-2\delta}\\
&~~~~~~~~~+\E\int_0^{t}|f_{h_l}(X^{\vv}_{h_l}(\eta_{h_l}(s)), X^{\vv}_{h_l}(\eta_{h_l}(s-\tau)))-f(X^{\vv}_{h_l}(\eta_{h_l}(s)), X^{\vv}_{h_l}(\eta_{h_l}(s-\tau)))|^2{\mbox d}s\\
&\le C\int_0^{t}\sup_{0\le u\le s}\E|X^\vv_{h_{l}}(u)-X_{h_{l-1}}^\vv(u)|^2{\mbox d}s+Ch_{l-1}^{2\delta}.
\end{split}
\end{equation*}
Since we have $\delta<1/2$, similar to the estimation of $J_1(t)$, we obtain
\begin{equation*}
\begin{split}
J_2(t)\le C\int_0^{t}\sup_{0\le u\le s}\E|X^\vv_{h_{l}}(u)-X_{h_{l-1}}^\vv(u)|^2{\mbox d}s+Ch_{l-2}^{2\delta}.
\end{split}
\end{equation*}
By \eqref{delones}, we compute
\begin{equation*}
\begin{split}
J_3(t)\le C\int_0^{t}\sup_{0\le u\le s}\E|X^{\vv}_{h_l}(u)-X^{\vv}_{h_{l-1}}(u)|^2{\mbox d}s.
\end{split}
\end{equation*}
By \eqref{delbound} and Lemma \ref{1gaosi}, we derive
\begin{equation*}
\begin{split}
J_4(t)\le Ch_{l}^{1-2\delta}+C\vv h_{l}^{\frac{1-2\delta}{2}}.
\end{split}
\end{equation*}
With \eqref{delbound} again,
\begin{equation*}
\begin{split}
J_5(t)\le Ch_{l}^{1-2\delta}.
\end{split}
\end{equation*}
By \eqref{dec}, Lemma \ref{1pmoment} and Lemma \ref{1miss}, we derive
\begin{equation*}
\begin{split}
J_6(t)\le C\int_0^{t}\sup_{0\le u\le s}\E|X^{\vv}_{h_l}(u)-X^{\vv}_{h_{l-1}}(u)|^2{\mbox d}s+Ch_{l-1}^{2-2\delta}+\vv^2h_{l-1}.
\end{split}
\end{equation*}
By sorting those inequalities and using the Gronwall inequality, the desired result can be obtained.\hfill $\Box$

\begin{remark}
{\rm Under one-sided Lipschitz condition, we show that the second moment of two coupled paths is bounded by $\mathcal{O}(h_{l-1}^{1-2\delta}+h_{l-2}^{2\delta}+\vv^2h_{l-1})$. By Theorem \ref{1election}, as the procedure of Theorem \ref{th1}, we can show that under assumptions (H1) and (H2), we have
\begin{equation*}
\begin{split}
\sup_{0\le n<M^{l-1}}{\rm Var}(\Psi(X_{h_l}^\vv(t_n))-\Psi(X_{h_{l-1}}^\vv(t_n)))\le Ch_{l-1}^{1-2\delta}+Ch_{l-2}^{2\delta}+C\vv^2h_{l-1}.
\end{split}
\end{equation*}
Different from the global Lipschitz case, the efficiency can not be improved compared with the second moment since the drift coefficient is one-sided Lipschitz not global Lipschitz. }
\end{remark}

\end{document}